\documentclass[11pt,a4paper]{article}
\usepackage{fontspec}
\newcommand{\latinfont}{\rmfamily}
\usepackage{amsmath,amssymb,amsthm,mathtools}
\usepackage[margin=24mm,headheight=15pt]{geometry}
\usepackage[unicode,hidelinks]{hyperref}
\usepackage{fancyhdr}
\newtheoremstyle{volume}{7pt}{7pt}{\normalfont}{}{\bfseries}{.}{.5em}{}
\theoremstyle{volume}
\newtheorem{theorem}{Theorem}[section]
\newtheorem{lemma}[theorem]{Lemma}
\newtheorem{proposition}[theorem]{Proposition}
\newtheorem{corollary}[theorem]{Corollary}

\newtheorem{example}[theorem]{Example}

\newcommand{\R}{\mathbb R}
\newcommand{\Q}{\mathbb Q}
\newcommand{\Z}{\mathbb Z}
\newcommand{\cG}{\mathcal{G}}
\newcommand{\cH}{\mathcal{G}_{\mathrm{ex}}}
\newcommand{\cL}{\mathcal{L}}
\newcommand{\totalVol}{\mathrm{Vol}}
\DeclareMathOperator{\Diff}{Diff}
\DeclareMathOperator{\Flux}{Flux}
\DeclareMathOperator{\SO}{SO}
\DeclareMathOperator{\SL}{SL}
\DeclareMathOperator{\Map}{Map}

\DeclareMathOperator{\Vol}{Vol}
\DeclareMathOperator{\ev}{ev}
\DeclareMathOperator{\pr}{pr}
\DeclareMathOperator{\im}{im}
\DeclareMathOperator{\Hom}{Hom}
\DeclareMathOperator{\Sym}{Sym}

\title{\textbf{Integral periods of the top volume class on volume-preserving diffeomorphism groups}}
\author{Takefumi Nosaka\\[2pt]\small Department of Mathematics, Institute of Science Tokyo}\date{November 18, 2026}
\hypersetup{pdfauthor={Takefumi Nosaka},pdftitle={Integral periods of the top volume class on volume-preserving diffeomorphism groups}}
\begin{document}
\maketitle\thispagestyle{plain}
\begin{abstract}
Let $M$ be a closed connected oriented aspherical $n$-manifold ($n\ge3$) with centerless
fundamental group, vanishing positive-degree rational Pontryagin classes and positive
volume form $v$. We normalize the periods of the top volume class on the discrete group
$\operatorname{Diff}_{v,0}(M)$ and its flux kernel by $\int_M v$. For odd $n$, $\mathbb{Q}$ lifts
additively to both integral homology groups. For even $n$, Betti numbers determine each
period group, either $\{0\}$ or $\mathbb{R}$. In dimension three the flux kernel's periods form
a real closed field, with additive lifts on a real closed subfield. For hyperbolic $M^3$
with hyperbolic volume form $v$, the Brandenbursky--Marcinkowski bounded volume class
on $\operatorname{Homeo}_0(M,v)$ is nonzero.
\end{abstract}
\begingroup
\par
\endgroup

\setcounter{section}{-1}
\section{Introduction}\label{sec:intro}

When the holonomy of a flat bundle preserves a volume form on the fiber, the form defines a transverse
volume class on the total space. Here this class has a component in the cohomology of the base, which
we call the top volume class. Which real numbers occur as its periods on integral homology? Can one
choose an integral class for each period so that the choices respect addition? We consider these
questions for the identity component of the volume-preserving diffeomorphism group and its flux
kernel.

Let $M^n$ be a closed connected oriented smooth aspherical manifold with $Z(\pi_1M)=1$.
Here $Z$ denotes the center of a group, and asphericity means that the universal cover is contractible.
Fix a smooth volume form $v$ compatible with the orientation, and write $\totalVol\coloneqq\int_Mv>0$ for the total volume.
We use ordinary composition and the $C^\infty$ topology on the group of $v$-preserving diffeomorphisms.
Let $\cG$ denote its identity component and $\cH$ the kernel of the volume flux homomorphism on $\cG$, so that
\[
 \cG\coloneqq\Diff_{v,0}(M),\qquad
 \cH\coloneqq\ker\bigl(\Flux:\cG\to H^{n-1}(M;\R)\bigr).
\]
Since the fundamental group is centerless, the flux group vanishes, so $\Flux$ is a continuous
homomorphism on $\cG$ \cite[Theorem~15(2)]{KKM}. Its surjectivity (Section~\ref{subsec:flat}) and Poincar\'e duality give $\cG/\cH\cong H_1(M;\R)$. Thus $\cG=\cH$ if and
only if $H_1(M;\R)=0$. For a group $H$, write $H^\delta$ for the group with the discrete topology and $BH^\delta$ for its
classifying space. Write $b_j(M)\coloneqq\dim_\Q H^j(M;\Q)$ for the $j$-th rational Betti number.

The universal flat $M$-bundle $\pi:\mathcal M_{\cG}\to B\cG^\delta$ carries a transverse volume class
$\Theta\in H^n(\mathcal M_{\cG};\R)$ whose restriction to each fiber is $[v]$.
In Section~\ref{subsec:flat} we construct a map $r:\mathcal M_{\cG}\to M$ such that $(\pi,r)$ is a fiber homotopy trivialization.
The \emph{top volume class} is the component of $\Theta$ of bidegree $(n,0)$ in the resulting K\"unneth decomposition.

Let $\pi_!$ denote integration along the fiber, using the base-first convention for the product orientation.
The top volume class and its integral period homomorphism are given by
\begin{equation}\label{intro:xi}
 \xi_{M,n}\coloneqq \totalVol^{-1}\pi_!(\Theta\smile r^*[v]),\qquad
 P_{\cL}:H_n(B\cL^\delta;\Z)\longrightarrow\R,\quad
 z\longmapsto\langle\xi_{M,n}|_{\cL},z\rangle
 \quad(\cL=\cG,\cH).
\end{equation}
We call $\im P_{\cL}$ the \emph{integral period group} and $P_{\cL}(z)/\totalVol$ a \emph{normalized integral period}. By the
universal coefficient theorem over the reals, $\xi_{M,n}|_{\cL}$ vanishes if and only if $P_{\cL}=0$. Throughout, $i:\cH\hookrightarrow \cG$ is the
inclusion and $z_{\R}$ is the image of the integral homology class $z$ in real homology. The $\ell^1$-seminorm $\|\cdot\|_1$
of a real homology class is the infimum, over all real singular cycles representing it, of the sum of
the absolute values of their coefficients. Write $\|M\|=\|[M]_\R\|_1$ for the simplicial volume.

A \emph{coherent family of periods} indexed by $\Q$ is a family of integral homology classes $(z_t)_{t\in\Q}$
such that $P_{\cL}(z_t)=\totalVol\,t$ and $z_{s+t}=z_s+z_t$ for all $s,t\in\Q$.

McDuff constructed canonical classes and established the comparison with the section space of Haefliger structures
with transverse volume \cite{McDuffI,McDuffIII,McDuffCan}. For a topological group $H$, we call the homotopy fiber
of $BH^\delta\to BH$ its \emph{local homology space}. McDuff studied the nonvanishing of the top-degree canonical class
on this space for groups of volume-preserving diffeomorphisms of odd-dimensional spheres \cite[p.~346]{McDuffCan}.
Moskowitz treated parallelizable manifolds and odd-dimensional stably parallelizable manifolds \cite{Moskowitz}.
Hurder constructed measure-theoretic characteristic classes for
Haefliger structures with transverse volume and proved that, in the framed case, simultaneous
evaluation of specified characteristic classes on integral homology is surjective \cite[Theorems~5.1--5.2, Corollary~5.3]{Hurder}. Here we use
McDuff's integral homology comparison to express evaluation of the top volume class in terms of point
evaluation on the section space and to determine periods and their integral lifts.

We say that $M$ is \emph{rationally Pontryagin-trivial} if $p_i(TM)=0\in H^{4i}(M;\Q)$ for every $i\ge1$.
In odd dimensions this hypothesis yields a coherent family indexed by $\Q$ on the flux kernel itself.

\begin{theorem}[Coherent families of integral periods in odd dimensions]\label{thm:odd}
Let $n\ge3$ be odd, and let $M^n$ be a closed connected oriented aspherical manifold with centerless fundamental group.
If $M$ is rationally Pontryagin-trivial, then for every positive volume form $v$ there is an additive injection
$\alpha_{\cH}^{\Q}:\Q\to H_n(B\cH^\delta;\Z)$. The composite $\alpha_{\cG}^{\Q}\coloneqq i_*\alpha_{\cH}^{\Q}$ is also injective, and
\[
 P_{\cH}(\alpha_{\cH}^{\Q}(t))=P_{\cG}(\alpha_{\cG}^{\Q}(t))=\totalVol\,t
 \qquad(t\in\Q).
\]
The images of the two maps are direct summands of their respective integral homology groups.
Moreover, $\|(\alpha_{\cH}^{\Q}(t))_\R\|_1$ and $\|(\alpha_{\cG}^{\Q}(t))_\R\|_1$ are both at least $\|M\|\,|t|$.
\end{theorem}

In even dimensions, expanding the square-zero relation for the transverse volume class in its K\"unneth components
gives a quadratic identity relating the top volume class to the mixed terms. The component of bidegree $(1,n-1)$
represents flux and hence vanishes on the flux kernel. This accounts for the two different ranges of Betti numbers
in the following criterion.

\begin{theorem}[The integral period criterion in even dimensions]\label{thm:even}
Let $M$ be a closed connected oriented smooth aspherical manifold of dimension $n=2m\ge4$ with centerless fundamental group.
For $1\le a\le m$, put $b_{[a,m]}(M)\coloneqq\sum_{p=a}^{m}b_p(M)$.
Fix a positive volume form $v$ and take $(\cL,a)$ to be either $(\cG,1)$ or $(\cH,2)$.
If $b_{[a,m]}(M)=0$, then $\xi_{M,n}|_{\cL}=0$ and $\im P_{\cL}=\{0\}$.
If $b_{[a,m]}(M)>0$ and $M$ is rationally Pontryagin-trivial, then $\im P_{\cL}=\R$.
\end{theorem}

If $M$ is rationally Pontryagin-trivial, rational cohomology in intermediate degrees also gives a
sufficient condition for every real number to occur as a period in odd dimensions (Proposition~\ref{oddreal:main}).
If the manifold admits a nonpositively curved metric and has a nonzero rational Pontryagin class of
positive degree, the top volume class vanishes on both groups (Proposition~\ref{scope:bott}).

In dimension three, the real closed fields $\mathcal P_{\mathrm{div}}\subseteq\mathcal P\subseteq\R$
constructed in Section~\ref{sec:threeD-refinement} are independent of $M$ and $v$.
Proposition~\ref{h3:period-characterization} identifies $\mathcal P$ with the normalized integral periods of the flux kernel.
A normalized period $u$ belongs to $\mathcal P_{\mathrm{div}}$ precisely when there is a homomorphism
$r\mapsto z_r$ from $\Q$ to $H_3(B\cH^\delta;\Z)$ such that $z_r$ has normalized period $ru$ for every
$r\in\Q$.
Integral lifts can be chosen additively on all of $\mathcal P_{\mathrm{div}}$
(Theorem~\ref{h3:geo:integer-additive-section}).
Let $E_{\log}$ be the smallest real closed subfield of $\R$ containing $\log b/\log a$ for all real algebraic $a,b>1$.
It is countable, being the real closure in $\R$ of the countable field generated by these ratios.
A construction using the Sol group gives $E_{\log}\subseteq\mathcal P_{\mathrm{div}}$ (Theorem~\ref{h3:src:logfield}).

\begin{theorem}[Coherent families of periods in dimension three]\label{thm:threeD-E}
Let $M$ be a closed connected oriented aspherical three-manifold with centerless fundamental group,
and let $v$ be any positive volume form. There is a noncanonical additive injection $\alpha_{\cH}^{E_{\log}}:E_{\log}\to H_3(B\cH^\delta;\Z)$. The map $\alpha_{\cG}^{E_{\log}}\coloneqq i_*\alpha_{\cH}^{E_{\log}}$ is
also injective, and for $\cL=\cH,\cG$ we have
\[
 P_{\cL}(\alpha_{\cL}^{E_{\log}}(u))=\totalVol\,u,\qquad
 \|(\alpha_{\cL}^{E_{\log}}(u))_{\R}\|_1\ge\|M\|\,|u|
 \qquad(u\in E_{\log}).
\]
\end{theorem}

Restricting this family to $\Q$ gives a choice of the family in Theorem~\ref{thm:odd}.

Brandenbursky and Marcinkowski asked whether their bounded volume class on the volume-preserving
homeomorphism group is nonzero for closed hyperbolic three-manifolds and for higher-dimensional
hyperbolic manifolds \cite[discussion after Theorem~1.1, p.~35]{BM}. We recall their class and specify the convention
for group multiplication in Proposition~\ref{common:bm-comparison}.

\begin{theorem}[The bounded volume class of a closed real hyperbolic manifold]\label{thm:bm-detection}
Let $M^n$ be a closed connected oriented real hyperbolic manifold of dimension $n\ge3$, and let $v$ be its
hyperbolic volume form. If $4\mid n$, assume that $M$ is not a rational homology sphere. Let $H=\operatorname{Homeo}_0(M,v)$ denote the
group of homeomorphisms joined to the identity by a $v$-preserving isotopy. Then the
Brandenbursky--Marcinkowski bounded volume class on $H$ is nonzero. In odd dimensions $H_n(BH^\delta;\Z)$ contains a direct
summand isomorphic to $(\Q,+)$ on which evaluation of the class is injective. In even dimensions evaluation on
$H_n(BH^\delta;\Z)$ is surjective onto $\R$.
\end{theorem}

In Section~\ref{sec:common} we construct the comparison with the section space and obtain a group cocycle
by averaging with respect to the volume form. For an integral class $z$ with normalized period $t$,
this cocycle gives the lower bound $\|z_\R\|_1\ge\|M\|\,|t|$.
For measure-theoretic induction in group cohomology, see Reznikov \cite[preprint, \S\S7.2--7.3]{Reznikov} and Brandenbursky--Marcinkowski
\cite{BMInd,BM}. Section~\ref{sec:models} develops the section-space models and the lifting obstruction used in Section~\ref{sec:odd} to
prove the odd-dimensional period theorem and the even-dimensional criterion.

Section~\ref{subsec:bm-proof} uses integral periods to prove Theorem~\ref{thm:bm-detection}. Section~\ref{sec:threeD-refinement} describes the normalized integral
periods of the flux kernel in dimension three and constructs a coherent family over $E_{\log}$ (Theorem~\ref{thm:threeD-E}).
Section~\ref{sec:applications} treats closed real hyperbolic manifolds with arbitrary positive volume forms.
When the outer action on the fundamental group has finite image, we use transfer to detect periods on
the full group of volume-preserving diffeomorphisms.

In Section~\ref{sec:scope} we apply Theorem~\ref{thm:even} to a four-dimensional torus bundle and
Proposition~\ref{oddreal:main} to a family of higher-rank locally symmetric spaces. We also give examples of
vanishing on closed complex hyperbolic manifolds. In the torus bundle, every real period on the identity
component is realized by an integral class whose image in real homology has zero $\ell^1$-seminorm.
Consequently, the top volume class has no bounded representative.
For closed real hyperbolic three-manifolds, Appendix~\ref{app:local} constructs a globally smooth bounded
cocycle with the same germ as a specified local point-evaluation cocycle. Its ordinary cohomology class
is the top volume class. Combined with the transfer calculation in Section~\ref{sec:applications}, this proves
that the transferred class considered by Nosaka \cite[\S6.1]{Nosaka} is nonzero.

\noindent\textbf{Notation and standing assumptions.} Unless otherwise stated, $M^n$ is a closed
connected oriented smooth aspherical manifold whose fundamental group $\Gamma\coloneqq\pi_1M$ has trivial center. Write $q_M:\widetilde M\to M$
for the universal covering map, with $\widetilde M$ contractible. We use $v,\totalVol,\cG,\cH$ as defined above.

\noindent\textbf{Acknowledgments.} The author benefited from discussions with OpenAI's ChatGPT while working on the arguments
concerning bounded cocycles and the proofs in the appendix, and the revise of the exposition.
The author takes full responsibility for the contents of this paper.

\section{The top volume class and the comparison map}\label{sec:common}

We extract the top volume class from the transverse volume class of the universal flat bundle and
use McDuff's comparison to identify its pullback to the homotopy fiber with the pullback of the
point-evaluation class from the section space.
Averaging cocycles against volume then gives a bound on periods. For real hyperbolic manifolds,
volumes of geodesic simplices give a bounded representative, which we compare with the bounded volume
class on the volume-preserving homeomorphism group.

\subsection{Flat bundles and comparison with the section space}\label{subsec:flat}
Suppose $n\ge3$. A volume-preserving isotopy generated by a time-dependent vector field $X_t$ has flux
$[\int_0^1\iota_{X_t}v\,dt]$. With the $C^\infty$ subspace topology, the kernel of the flux homomorphism,
$\cH=\Diff_{\mathrm{ex}}(M,v)$, is a Fr\'echet--Lie subgroup \cite[\S2, p.~345, items~1--2]{HTV}.

Let $\cL$ be a subgroup of $\Diff_v(M)$ whose elements are homotopic to the identity as maps $M\to M$.
Both $\cG$ and $\cH$ satisfy this condition. For $g\in\cL$, lift a homotopy from the identity to $g$,
starting at the identity of $\widetilde M$. The resulting lift $\widetilde g$ commutes with every deck transformation.
Any two such lifts differ by a central deck transformation. Since $\Gamma$ has trivial center, the lift is unique
and the assignment $g\mapsto\widetilde g$ respects composition.
The group $\cG$ is locally path connected, and evaluation sends $\pi_1\cG$ into the center of $\Gamma$.
Thus $(g,\widetilde x)\mapsto g(q_M(\widetilde x))$ admits a continuous lift on $\cG\times\widetilde M$.
The lift that restricts to the identity on $\{1\}\times\widetilde M$ sends $(g,\widetilde x)$ to $\widetilde g\widetilde x$.
Locally it is the composite of smooth evaluation with a smooth local inverse of $q_M$.
The lifted action is therefore smooth for the Fr\'echet--Lie group structure on $\cG$.

Let $\pi_{\cL}:\mathcal M_{\cL}\coloneqq E\cL^\delta\times_{\cL^\delta}M\to B\cL^\delta$ be the universal flat bundle. In the homogeneous bar model the forms $\pr_M^*v$ on the spaces $\Delta^k\times M$ are
compatible with the face maps and the group action, and their squares vanish. For each finitely
generated subgroup $\mathcal K\le\cL$, the multiplicative simplicial de Rham comparison
\cite[Theorems~6.4, 6.10]{Dupont} gives a class $\Theta_{\mathcal K}\in H^n(\mathcal M_{\mathcal K};\R)$ with $\Theta_{\mathcal K}^2=0$.
These classes are compatible with restriction. The projection of the support of a singular chain in
$\mathcal M_{\cL}$ lies in a finite subcomplex of $B\cL^\delta$ \cite[Proposition~A.1]{Hatcher}.
The cell labels of that subcomplex lie in a finitely generated subgroup $\mathcal K\le\cL$.
Thus the chain is supported in the pullback $\mathcal M_{\mathcal K}$ over $B\mathcal K^\delta$.
Homology is therefore the direct limit
over these subgroups. By the universal coefficient theorem, the classes $\Theta_{\mathcal K}$ determine a unique
real class $\Theta_{\cL}$. The same argument gives $\Theta_{\cL}^2=0$, since $\Theta_{\mathcal K}^2=0$ for every $\mathcal K$.
The canonical lifts identify $E\cL^\delta\times\widetilde M$ with the contractible universal cover and give $\pi_1\mathcal M_{\cL}\cong\Gamma\times \cL^\delta$. Since $M$ is a
model for $B\Gamma$, there is a map $r_{\cL}:\mathcal M_{\cL}\to M$ inducing projection onto the first factor. Any map homotopic to the
identity on the reference fiber induces a homomorphism sending the second factor into the center of
$\Gamma$. Thus $r_{\cL}$ is unique up to homotopy under this condition, and $(\pi_{\cL},r_{\cL})$ gives a fiber homotopy
trivialization.

Since $M$ has finite CW type, we may write the K\"unneth decomposition in this trivialization as
\begin{equation}\label{common:kun}
 \Theta_{\cL}=\xi_{\cL}\times1+1\times[v]+\Theta_{\cL,\mathrm{mix}}.
\end{equation}
The last term is the sum of the components of positive degree in both base and fiber. The pure base
component $\xi_{\cL}$ equals $\totalVol^{-1}(\pi_{\cL})_!(\Theta_{\cL}\smile r_{\cL}^*[v])$ and is natural under subgroup inclusions. In particular $\xi_{\cG}=\xi_{M,n}$ and $\xi_{\cH}=\xi_{M,n}|_{\cH}$. For general
$\cL$ we define the integral period by $P_{\cL}(z)\coloneqq\langle\xi_{\cL},z\rangle$, as in \eqref{intro:xi}. Reversing the order of the factors gives
\begin{equation}\label{common:sign}
 \frac{1}{\totalVol}(\pi_{\cL})_!\bigl(r_{\cL}^*[v]\smile\Theta_{\cL}\bigr)=(-1)^n\xi_{\cL}.
\end{equation}

We next relate the top volume class to point evaluation on the section space. Let $E_n$ be the
classifying space of the \'etale groupoid of germs of smooth local diffeomorphisms preserving the
standard volume form on $\R^n$. The functor sending a germ $[f]_x:x\to f(x)$ to its derivative
$D_xf\in\SL_n(\R)$ induces a map on classifying spaces. We use the fibration model
\begin{equation}\label{common:haefliger}
 F_n\xrightarrow{i_F} E_n
 \xrightarrow{\nu}B\SL_n(\R)\simeq B\SO(n).
\end{equation}
Write $u\in H^n(E_n;\R)$ for the universal transverse volume class and $u_F\coloneqq i_F^*u\in H^n(F_n;\R)$ for its restriction to the fiber.
Haefliger's connectivity result and McDuff's homotopy group calculation give, for $n\ge3$,
\begin{equation}\label{common:localhomotopy}
 F_n\text{ is }(n-1)\text{-connected},\qquad
 \pi_nF_n\xrightarrow[\langle u_F,-\rangle]{\cong}(\R,+),\qquad
 \pi_{n+1}F_n=0.
\end{equation}
For the connectivity and homotopy group calculations, see \cite[p.~609 and Corollary~1.3 (p.~612)]{McDuffI}.
Hurder identifies the displayed isomorphism with evaluation against $u_F$ \cite[p.~382]{Hurder}.

Let $\tau_M:M\to B\SL_n(\R)$ classify the tangent bundle equipped with $v$.
An atlas in which $v$ is the standard volume form determines a section $s_0$ of $\tau_M^*E_n\to M$.
Let $\mathcal S_M$ be the component of the section space containing $s_0$.
The corresponding lift $\sigma_0:M\to E_n$ satisfies $\sigma_0^*u=[v]$.
Fix $a_0\in M$ and a volume frame at $a_0$. This frame identifies evaluation at $a_0$ with a map
$\ev_{a_0}:\mathcal S_M\to F_n$. When necessary, we use a CW model for the section space and take the
universal cover $\rho:\widetilde{\mathcal S}_M\to\mathcal S_M$ in that model.

For $\cL=\cG,\cH$, let $B_{\rm top}\cL$ denote the classifying space of the topological group, and put
\[
 \bar B\cL\coloneqq\operatorname{hofib}(B\cL^\delta\to B_{\rm top}\cL),
 \qquad p_{\cL}:\bar B\cL\to B\cL^\delta.
\]
McDuff's comparison theorem and its flux-kernel version give maps
\begin{equation}\label{common:mcduff}
 m_{\cG}:\bar B\cG\longrightarrow\mathcal S_M,
 \qquad
 m_{\cH}:\bar B\cH\longrightarrow\widetilde{\mathcal S}_M
\end{equation}
inducing integral homology isomorphisms \cite[Theorem~1, pp.~531--533]{McDuffIII}. Here $\pi_1\mathcal S_M$ is identified with $H^{n-1}(M;\R)$, and the kernel in McDuff's
formulation agrees with $\cH$ because the flux group vanishes. The map $\bar Bi:\bar B\cH\to\bar B\cG$ induced by the inclusion $i$
satisfies $\rho m_{\cH}\simeq m_{\cG}\bar Bi$. Passing from the full group used by McDuff to the identity component preserves the
homotopy fiber, as one sees by comparing the two classifying spaces over the classifying space of
their common component group.

We record a decomposition of $\cG$ as a space. Choose volume-preserving vector fields $X_1,\ldots,X_k$ such that the
classes $[\iota_{X_j}v]$ form a basis of $H^{n-1}(M;\R)$. For $a=\sum_{j=1}^k a_j[\iota_{X_j}v]$, define $s(a)$
by composing the time-$a_j$ flows of $X_j$ in the chosen order. Then $\Flux(s(a))=a$, and $s$ is continuous.
It follows that $(h,a)\mapsto h\,s(a)$ defines a homeomorphism $\cH\times H^{n-1}(M;\R)\to\cG$.

\begin{proposition}[The global comparison]\label{common:comparison}
For the identity component, the comparison gives
\begin{equation}\label{eq:comparison}
 p_{\cG}^*\xi_{M,n}=m_{\cG}^*\ev_{a_0}^*u_F.
\end{equation}
For the flux kernel, the corresponding identity is
$p_{\cH}^*\xi_{\cH}=m_{\cH}^*(\ev_{a_0}\rho)^*u_F$.
\end{proposition}
\begin{proof}
By the definition of $\bar B\cG$ as the homotopy fiber of $B\cG^\delta\to B_{\rm top}\cG$, the
pulled-back bundle is trivialized as $\bar B\cG\times M$. The centerless hypothesis and the condition on the reference fiber give $r\simeq\pr_M$ in
these coordinates. Restricting the pulled-back transverse volume class to $\bar B\cG\times\{a_0\}$ extracts its pure base
component $p_{\cG}^*\xi_{M,n}$. By McDuff's construction this restriction is $m_{\cG}^*\ev_{a_0}^*u_F$ \cite[Proof of Lemma~3.1]{McDuffIII}. The assertion for $\cH$ follows from
$\rho m_{\cH}\simeq m_{\cG}\bar Bi$.
\end{proof}

Write $\mathsf T_{\cL}\coloneqq p_{\cL*}(m_{\cL*})^{-1}$ for the map to $H_*(B\cL^\delta;\Z)$.
Its domain is $H_*(\mathcal S_M;\Z)$ for $\cL=\cG$ and $H_*(\widetilde{\mathcal S}_M;\Z)$ for $\cL=\cH$.
In either case, Proposition~\ref{common:comparison} identifies the period of $\mathsf T_{\cL}(z)$ with the pairing of $z$
with the corresponding point-evaluation class.

\subsection{Volume cocycles and simplicial volume}\label{subsec:average}
A homogeneous $d$-cochain is a function $c:\cL^{d+1}\to\R$ invariant under simultaneous left
translation. Omitting each argument in turn gives the alternating sum for its coboundary $\delta c$.
The corresponding inhomogeneous cochain is $f(g_1,\ldots,g_d)=c(1,g_1,g_1g_2,\ldots,g_1\cdots g_d)$.
Ordinary group cohomology requires no regularity of cochains. A cochain is locally smooth if it is
smooth near $(1,\ldots,1)$. Borel measurability and boundedness refer to the whole domain.

\begin{proposition}[A cocycle representative obtained by averaging]\label{common:continuous}
Let $b$ be an alternating homogeneous $n$-cocycle on $\Gamma$ representing $[v]$ under the identification $B\Gamma\simeq M$.
From $b$ one can construct a homogeneous cocycle $C_b$ representing $\xi_{\cL}$.
If $b$ is bounded, then $\|C_b\|_\infty\le\|b\|_\infty$. For $\cL=\cG$, the cocycle $C_b$ is globally smooth.
\end{proposition}
\begin{proof}
Let $\Delta(\Gamma)$ be the simplex of finitely supported probability vectors on $\Gamma$. Since $M$ is a
finite-dimensional model for $B\Gamma$, the group $\Gamma$ is torsion-free \cite[Proposition~2.45]{Hatcher}. If left translation by $\gamma\in\Gamma$ preserves the finite support $S$ of a probability vector,
then for any $s\in S$ there is a $k>0$ with $\gamma^k s=s$. Hence $\gamma^k=1$, and torsion-freeness gives $\gamma=1$.
Thus $\Delta(\Gamma)$ is a free $\Gamma$-complex. Linear contraction to a vertex shows that it is contractible. Extend $b$ multilinearly to finitely supported
vectors. For a probability vector $p$ and tangent vectors $w_j$ with coordinate sum zero, set $\omega_b|_p(w_1,\ldots,w_n)\coloneqq n!b(p,w_1,\ldots,w_n)$.
Since $\delta b=0$, this expression is independent of $p$ and defines a $\Gamma$-invariant constant form compatible
with all faces. Its integral over the affine simplex with vertices $p_0,\ldots,p_n$ is $b(p_0,\ldots,p_n)$. In particular,
integration along the affine map that sends bar vertices to their unit vectors gives $b$, so $\omega_b$
represents the class on the quotient $B\Gamma$ corresponding to $[v]$.

Choose a smooth nonnegative function of compact support on $\widetilde M$ such that at each point
at least one of its $\Gamma$-translates is positive. Dividing the translates by their sum gives an equivariant map
$p:\widetilde M\to\Delta(\Gamma)$. On a bar simplex $(g_0,\ldots,g_k)$ put $p_i\coloneqq p(\widetilde g_i^{-1}\widetilde x)$.
The affine combination $\sum_i t_i p_i$ is compatible with the face maps and the group action. It therefore defines
a map $R:\mathcal M_{\cL}\to B\Gamma$ whose restriction to the reference fiber is the standard equivalence.
By the uniqueness of $r_{\cL}$, we have $R^*[\omega_b]=r_{\cL}^*[v]$.

To compute the fiber integral in \eqref{common:sign}, consider more generally a map $R$ to a space
equipped with a closed $n$-form $\omega$. Suppose that $R^*\omega$ is a simplicial differential form on
the bar model, compatible with the face maps and representing $r_{\cL}^*[v]$. Then
\begin{equation}\label{common:fiber-average}
 (g_0,\ldots,g_n)\longmapsto
 \int_{\Delta^n\times M}R^*\omega\wedge\pr_M^*v
\end{equation}
is a homogeneous cocycle representing $(\pi_{\cL})_!(r_{\cL}^*[v]\smile\Theta_{\cL})$. Indeed, $\Theta_{\cL}$ is represented by $\pr_M^*v$ in a flat trivialization.
On restriction to any finitely generated subgroup, multiplicativity of the simplicial de Rham
comparison \cite[Theorems~6.4, 6.10]{Dupont} identifies the wedge of the two forms with their cup product in cohomology. Integrating
first over the fiber $M$ and then over $\Delta^n$ gives \eqref{common:fiber-average}. The direct limit description in Section~\ref{subsec:flat} gives
the same formula on the whole group. By \eqref{common:sign} and the affine simplex integral for $\omega=\omega_b$, we obtain
\begin{equation}\label{common:partitionaverage}
 C_b(g_0,\ldots,g_n)=\frac{(-1)^n}{\totalVol}\int_M b(p_0,\ldots,p_n)\,v(x).
\end{equation}
The integrand is invariant under deck transformations, and a volume-preserving change of variables
gives homogeneity. Each $p_i$ is a probability vector, which yields the claimed norm bound for bounded
$b$.

For $\cL=\cG$ the action by canonical lifts and the group inverse are globally smooth. Cover $M$ by finitely
many coordinate neighborhoods with relatively compact lifts and split the integral using a
subordinate smooth partition of unity. Near any fixed tuple of group elements, the points $\widetilde g_i^{-1}\widetilde x$ on
each chosen lift lie in a common compact set. By local finiteness of the partition defining $p$, only
finitely many coordinates occur there. Differentiation under the integral over a compact domain now shows
that $C_b$ is smooth, even when $b$ is unbounded.
\end{proof}

\begin{corollary}[A period bound from simplicial volume]\label{common:lone}
For $\cL$ as above and any $z\in H_n(B\cL^\delta;\Z)$, we have
\[
 \|z_{\R}\|_1\ge \frac{\|M\|}{\totalVol}\,|P_{\cL}(z)|.
\]
If $\|M\|>0$, the class $\xi_{\cL}$ has a bounded representative of supremum norm at most $\totalVol/\|M\|$.
\end{corollary}
\begin{proof}
The bound is immediate if $\|M\|=0$. Suppose $\|M\|>0$. By Frigerio \cite[Proof of Lemma~6.1]{FrigerioBC} there is a bounded singular cocycle taking
the value $\|M\|$ on $[M]_\R$ and having supremum norm at most $1$. Multiplying by $\totalVol/\|M\|$ gives a bounded singular
representative of $[v]$ since $H_n(M;\R)=\R[M]_\R$. Contractibility of $\widetilde M$ gives a $\Gamma$-equivariant map from the
homogeneous bar model $E\Gamma$ to $\widetilde M$, inducing the standard equivalence $B\Gamma\simeq M$ on quotients.
Evaluate the singular cocycle on the images of the bar simplices and take the signed average over all
permutations of the vertices. This gives an alternating group cocycle $b$ representing $[v]$ with the same norm bound. Proposition~\ref{common:continuous} gives a representative
of $\xi_{\cL}$ satisfying $\|C_b\|_\infty\le\totalVol/\|M\|$.
Pulling $C_b$ back along Frigerio's vertex map \cite[Lemma~5.2]{FrigerioBC} gives a singular cocycle
of norm at most $\|C_b\|_\infty$. The underlying
chain map is an augmentation-preserving comparison of free resolutions, so the resulting singular
cocycle represents $\xi_{\cL}$ in ordinary cohomology. Hence $|P_{\cL}(z)|\le\|C_b\|_\infty\|z_{\R}\|_1$. If the manifold admits a hyperbolic metric,
then $\|M\|=V_{\rm hyp}/v_n$, where $V_{\rm hyp}$ is the volume of the metric of curvature $-1$ and $v_n$ is the supremum of the
volumes of hyperbolic $n$-simplices \cite[Theorem~7.4]{FrigerioBC}.
\end{proof}

Let $M=\Gamma\backslash\mathbb H^n$ be a closed real hyperbolic manifold, with $v$ the hyperbolic volume form for curvature $-1$. In
Proposition~\ref{common:hyperbolic} we take the map $R$ in \eqref{common:fiber-average} to be given by geodesic simplices and put $\omega=v$. For each fixed
$x\in M$, the integral of $R(\,\cdot\,,x)^*v$ over $\Delta^n$ is the signed volume of the corresponding geodesic simplex in $\mathbb H^n$
\cite[\S8.7]{FrigerioBC}. Integrating this quantity in $x$ against the volume measure $v$ gives \eqref{common:average}. We write $\Vol_{\mathbb H^n}[z_0,\ldots,z_n]$ for the
signed volume of the geodesic simplex with vertices $z_0,\ldots,z_n$.

\begin{proposition}\label{common:hyperbolic}
For $g_0,\ldots,g_n\in\cG$, define
\begin{align}
 A_{M,n}(g_0,\ldots,g_n)
 &\coloneqq\int_M\Vol_{\mathbb H^n}
 [\widetilde g_0^{-1}\widetilde x,\ldots,
   \widetilde g_n^{-1}\widetilde x]\,v(x),\label{common:average}\\
 F_{M,n}&\coloneqq\frac{(-1)^n}{\totalVol}A_{M,n}.\label{common:normalized}
\end{align}
Then $F_{M,n}$ is a globally smooth homogeneous group cocycle, and
\begin{equation}\label{common:bounded}
 [F_{M,n}]=\xi_{M,n},\qquad \|F_{M,n}\|_\infty\le v_n.
\end{equation}
\end{proposition}
\begin{proof}
In the hyperboloid model, for future-pointing unit vectors $z_i$ and $t\in\Delta^k$, set
\[
 w=\sum_i t_i z_i,\qquad
 \operatorname{st}(t;z_0,\ldots,z_k)=\frac{w}{\sqrt{-\langle w,w\rangle}}.
\]
The inequality $-\langle w,w\rangle=\sum_{i,j}t_it_j\cosh d(z_i,z_j)\ge(\sum_i t_i)^2=1$ shows that straightening is jointly smooth in $t$ and the vertices, even for
degenerate simplices. It is isometry-equivariant and respects faces. Thus $R(t,x)\coloneqq q_M(\operatorname{st}(t;\widetilde g_0^{-1}\widetilde x,\ldots,\widetilde g_k^{-1}\widetilde x))$ is compatible with faces
and the diagonal left action and descends to $\mathcal M_{\cG}$. It is the identity on the reference fiber, hence
homotopic to $r$. Putting $\omega=v$ in \eqref{common:fiber-average} gives $A_{M,n}$, a homogeneous representative of $\pi_!(r^*[v]\smile\Theta)=(-1)^n\totalVol\,\xi_{M,n}$. The absolute volume
of each simplex is at most $v_n$, so $\|A_{M,n}\|_\infty\le\totalVol\,v_n$. As in Proposition~\ref{common:continuous}, smoothness of the lifted action and of
inversion implies that integration over the compact domain $\Delta^n\times M$ depends smoothly on the group
variables.
\end{proof}

\begin{proposition}[Comparison with the bounded volume class]\label{common:bm-comparison}
Equip $H=\operatorname{Homeo}_0(M,v)$ with the product $f\star h=h\circ f$, and define the homomorphism $J:\cG\to H$ by $J(g)=g^{-1}$. Let $\mathcal B_M\in H_b^n(H^\delta;\R)$ be the class induced from
the bounded hyperbolic volume class of the fundamental group using the measure $v$ \cite[\S\S2.2, 3.1--3.2]{BM}. Writing $\operatorname{comp}:H_b^n(H^\delta;\R)\to H^n(H^\delta;\R)$ for
the comparison map to ordinary cohomology, we have
\begin{equation}\label{common:bm}
 J^*\mathcal B_M=[A_{M,n}]_b,\qquad
 J^*\operatorname{comp}(\mathcal B_M)=(-1)^n\totalVol\,\xi_{M,n}.
\end{equation}
\end{proposition}
\begin{proof}
Choose a Borel section $S:M\to\mathbb H^n$ and a lift $z_*\in\mathbb H^n$ of the basepoint. Define $\gamma(f,x)\in\Gamma$ by $\widetilde fS(x)=\gamma(f,x)S(fx)$. A representative $C$ of $J^*\mathcal B_M$
is obtained by integrating the signed volume of the geodesic simplex with vertices $p_i(x)=\gamma(g_i^{-1},x)z_*$ against $v(x)$ \cite[Example~3.2, \S3.2]{BM}.
The cocycle $A_{M,n}$ in Proposition~\ref{common:hyperbolic} is given by the same integral with vertices $q_i(x)=\widetilde g_i^{-1}S(x)$. Using the
standard prism between the two vertex lists, put
\begin{equation}\label{common:bm-prism}
 \Pi(g_0,\ldots,g_{n-1})
 \coloneqq \int_M \sum_{j=0}^{n-1}(-1)^j
 \Vol_{\mathbb H^n}[p_0(x),\ldots,p_j(x),q_j(x),\ldots,q_{n-1}(x)]\,v(x).
\end{equation}
The integrand is Borel in $x$. It is an alternating sum of the signed volumes of $n$ geodesic
$n$-simplices, so its absolute value is at most $n v_n$. Hence $\|\Pi\|_\infty\le n v_n\totalVol$. Simultaneous left translation $g_i\mapsto hg_i$ applies
the same deck transformation $\gamma(h^{-1},x)$ to both lists of vertices at $h^{-1}x$. Invariance of $v$ therefore makes $\Pi$
homogeneous. Thus $\Pi\in C_b^{n-1}(\cG^\delta;\R)$. Since signed geodesic volume is a cocycle, the prism boundary formula gives
$\delta\Pi=A_{M,n}-C$ in this bounded cochain complex. This proves $[A_{M,n}]_b=[C]_b=J^*\mathcal B_M$, and Proposition~\ref{common:hyperbolic} gives the equality of
ordinary classes.
\end{proof}

If one uses the probability measure $v/\totalVol$ in the induction defining $\mathcal B_M$, both right-hand sides
of \eqref{common:bm} acquire the factor $1/\totalVol$.

\section{Section-space models and lifting obstructions}\label{sec:models}

We first rationalize the fibers of the bundle defining $\mathcal S_M$ and describe how to transport
the reference section along a homotopy of the classifying map of $TM$. We then give a homological
criterion for the vanishing of the degree-$2n$ lifting obstruction and use it to realize individual
real periods. In dimension three we also use the underlying relative Hurewicz argument with integral
coefficients.

\subsection{Fiberwise rationalization and the reference section}\label{subsec:section-models}

Assume $n\ge3$. We use CW models for classifying spaces and homotopy fibers and give mapping and section
spaces the compactly generated compact-open topology. For a nilpotent space $X$, write $\ell_X:X\to X_\Q$ for its
rationalization, abbreviated to $\ell$ when the space is clear. Write $K(A,j)$ for a connected Eilenberg--Mac
Lane space with homotopy group $A$ in degree $j$ and all other positive-degree homotopy groups zero.
Set $K_{\Q,n}\coloneqq K(\Q,n)$.

We choose a fibration model for the map whose components are the rationalized derivative map and the
classifying map of the transverse volume class. We write this model as
\begin{equation}\label{common:joint}
 (\nu_\Q,u_\Q):(E_n)_\Q\longrightarrow B\SO(n)_\Q\times K(\R,n).
\end{equation}
We also write $u_\Q$ for the real cohomology class corresponding to $u$. In the path-space model the projection
sends a path to its endpoint. To define $\ell_{E_n}$ in this model, use the chosen homotopy in the rationalization
square for the derivative coordinate and the constant path for the volume coordinate.
Then $\nu_\Q\ell_{E_n}=\ell_{B\SO(n)}\nu$ holds strictly. The composite of \eqref{common:joint} with projection
onto the first factor is the fibration $\nu_\Q$. By rationalization of \eqref{common:haefliger}, its fiber has
the homotopy type of $(F_n)_\Q$ \cite[Theorem~6.2.6]{MayPonto}.
The pullback $(\ell\tau_M)^*((E_n)_\Q)\to M$ is the fiberwise rationalization of $\tau_M^*E_n\to M$.
Let $\mathcal T_M$ be the component of its section space containing the section defined by
$\sigma_{0,\Q}\coloneqq\ell_{E_n}\sigma_0$. Using the same point and volume frame as in Section~\ref{subsec:flat},
we obtain $\ev_{a_0}:\mathcal T_M\to(F_n)_\Q$. Here $(u_F)_\Q$ denotes the real cohomology class
on $(F_n)_\Q$ that pulls back to $u_F$ under rationalization. Its pullback by $\ev_{a_0}$ is the point-evaluation class.

\begin{lemma}[Changing the tangent coordinate relative to the reference section]\label{common:tangent-transport}
Let $\mathcal H$ be a homotopy from $\ell\tau_M$ to $\tau':M\to B\SO(n)_\Q$. Lift $\mathcal H$ starting at $\sigma_{0,\Q}$ with the volume coordinate fixed, and
denote the terminal lift by $\sigma'_0$. Let $Y$ be a CW complex based at a $0$-cell, and let $\pr_M:Y\times M\to M$
be the projection. Suppose we are given a lift of $\tau'\circ\pr_M$ that agrees with $\sigma'_0$ on $\{*\}\times M$.
Then $\ell\tau_M\circ\pr_M$ has a lift with the same volume coordinate that agrees with $\sigma_{0,\Q}$ on $\{*\}\times M$.
\end{lemma}
\begin{proof}
Lift $\pr_M^*\mathcal H^{-1}$ along \eqref{common:joint} by the relative homotopy lifting property \cite[Lemma~1.3.2]{MayPonto}, keeping the volume coordinate fixed
and prescribing the reverse of the original lift on $\{*\}\times M$.
\end{proof}

Since $M$ has finite CW type and $F_n$ is simply connected, M{\o}ller's theorem identifies $\mathcal S_M\to\mathcal T_M$ as a
rationalization \cite[Theorems~4.1, 5.3]{Moller}. It commutes with point evaluation using the same frame and pulls $\ev_{a_0}^*(u_F)_\Q$ back to $\ev_{a_0}^*u_F$.
We take universal covers of CW models and work with based homotopy classes of families of sections in
these models. The induced map $\widetilde{\mathcal S}_M\to\widetilde{\mathcal T}_M$ rationalizes the higher homotopy groups and is therefore a rational
homology equivalence.

\medskip
We record the rational coordinates of the tangent bundle classifying map. In odd rank we have
\begin{equation}\label{odd:eq:bso}
 B\SO(2r+1)_\Q\simeq\prod_{j=1}^{r}K(\Q,4j)_{p_j}\text{.}
\end{equation}
In even rank $n=2m$ we have
\begin{equation}\label{even:bso}
 B\SO(2m)_\Q\simeq K(\Q,n)_e\times\prod_{j=1}^{m-1}K(\Q,4j)_{p_j}\text{.}
\end{equation}
In each case the polynomial generators of rational cohomology define the maps to the indicated factors
\cite[Chapter IV, \S5, Theorem~5.1]{BCM}. In even rank the top Pontryagin class satisfies $p_m=e^2$.
When $n=4k$, the classes $p_k$ and $e$ have the same degree but give distinct coordinates.
Thus, if $M$ is rationally Pontryagin-trivial, $\ell\tau_M$ is nullhomotopic in odd dimensions.
In even dimensions it is homotopic to the composite of the classifying map for $e(TM)$ with the inclusion
$j_e:K_{\Q,n}\to B\SO(n)_\Q$ of the Euler factor.

For $c>0$ define a dilation of the transverse coordinate $x\in\R^n$ by $d_c(x)=c^{1/n}x$. The map of germs $[f]_x\mapsto[d_c f d_c^{-1}]_{d_c(x)}$ preserves the
derivative, and $d_c$ pulls the standard volume form back to $c$ times itself. It therefore induces
self-maps $h_c$ of $E_n$ and $F_n$ satisfying
\begin{equation}\label{common:dilation}
 h_c^*u=cu,\qquad h_c^*u_F=cu_F.
\end{equation}
The map on $E_n$ commutes with $\nu$, so after rationalization it also acts on homotopy pullbacks formed
along $\nu_\Q$.

\subsection{A homological criterion for lifting}\label{subsec:obstruction}

For a rational vector space $A$ and an integer $r\ge2$, use the topological abelian group model of $K(A,r)$
obtained by geometric realization of a simplicial abelian group. Its integral homology groups in positive degrees
are $\Q$-vector spaces. With the Pontryagin product induced by addition, its rational homology algebra
is generated by $A$ in degree $r$ and takes the form
\begin{equation}\label{common:em-homology}
 H_*(K(A,r);\Q)\cong
 \begin{cases}
  \Sym_\Q A & (r\text{ even}),\\
  \bigwedge_\Q A & (r\text{ odd})
 \end{cases}.
\end{equation}
The finite-dimensional case follows from the calculation for $\dim_\Q A=1$ and the K\"unneth theorem. In
general, take the filtered colimit over the finite-dimensional subspaces in the simplicial
Eilenberg--Mac Lane model \cite[Theorem~5.2.8, Proposition~6.7.2]{MayPonto}.

We describe the relative obstruction in a form that will also be used in Section~\ref{sec:threeD-refinement}.
Let $p:D\to B$ be a map between simply connected spaces whose homotopy fiber is $(d-2)$-connected, where $d\ge3$.
The mapping cylinder pair $(M_p,D)$ is then $(d-1)$-connected. The homology exact sequence of this pair and
the relative Hurewicz theorem \cite[Theorem~4.32]{Hatcher} give the exact sequence
\begin{equation}\label{common:primary-obstruction}
 H_d(D;\Z)\xrightarrow{p_*} H_d(B;\Z)
 \longrightarrow H_d(M_p,D;\Z)\cong\pi_d(M_p,D).
\end{equation}
Let $(X,X_0)$ be a CW pair, and suppose that $f:X\to B$ has a prescribed lift on $X_0$.
Homotopy extension in the mapping cylinder gives a map of pairs $\bar f:(X,X_0)\to(M_p,D)$.
Since $(M_p,D)$ is $(d-1)$-connected, we can deform $\bar f$ relative to $X_0$ so that it sends the relative
$(d-1)$-skeleton into $D$. On each relative $d$-cell, the composite of the characteristic map with $\bar f$
represents the obstruction. Under the relative
Hurewicz isomorphism, evaluation of the primary obstruction is therefore $\bar f_*:H_d(X,X_0;\Z)\to H_d(M_p,D;\Z)$. If the term $\operatorname{Ext}^1_\Z(H_{d-1}(X,X_0;\Z),\pi_d(M_p,D))$ in the
universal coefficient sequence vanishes, this evaluation detects the vanishing of the obstruction
class. When $\dim X\le d$, the primary obstruction vanishes if and only if the prescribed lift extends in the
fibration model of $p$ \cite[\S4.3, pp.~416--419]{Hatcher}.

\begin{lemma}[The quadratic lifting criterion]\label{even:postnikov}
Suppose $n\ge3$, and let $Z$ be an $(n-1)$-connected rational space. Write $g:Z\to K(\pi_n Z,n)$ for the projection to its $n$-th
Postnikov stage. The lifting obstruction in degree $2n$ is given by a homomorphism
\begin{equation}\label{even:kclass}
 k:H_{2n}(K(\pi_n Z,n);\Q)\longrightarrow\pi_{2n-1}Z
\end{equation}
whose kernel is
\begin{equation}\label{even:kernelimage}
 \ker k=\im\bigl(g_*:H_{2n}(Z;\Q)\to H_{2n}(K(\pi_n Z,n);\Q)\bigr).
\end{equation}
The domain of $k$ is $\Sym^2_\Q(\pi_n Z)$ for even $n$ and $\bigwedge^2_\Q(\pi_n Z)$ for odd $n$.

Let $(X,X_0)$ be a CW pair with $\dim X\le2n$ and $\dim X_0\le n$, and let $f:X\to K(\pi_n Z,n)$ be a map. In the fibration model of $g$, a prescribed
lift of $f|_{X_0}$ extends over $X$ if and only if $f_*H_{2n}(X;\Q)\subseteq\im g_*$.
\end{lemma}

\begin{proof}
We proceed by induction from $P_nZ=K(\pi_n Z,n)$. Suppose $n<j\le2n-2$ and that the preceding stage $P_{j-1}Z$ splits as a product. By
\eqref{common:em-homology} the generators have degrees between $n$ and $j-1$, and products of two or more have degree at least
$2n$. Thus $H_{j+1}(P_{j-1}Z;\Q)=0$. Since $\pi_jZ$ is a rational vector space, the universal coefficient theorem implies that the
$k$-invariant at stage $j$ vanishes. The classification of Postnikov stages \cite[\S4.3]{Hatcher} therefore gives a
splitting
\begin{equation}\label{even:postnikovproduct}
 P_{2n-2}Z\simeq K(\pi_n Z,n)\times\prod_{j=n+1}^{2n-2}K(\pi_jZ,j)
\end{equation}
compatible with $g$. The same degree count shows that projection onto the first factor induces an
isomorphism on $H_{2n}(-;\Q)$.

The homotopy fiber of $p:Z\to P_{2n-2}Z$ is $(2n-2)$-connected, and $\pi_{2n}(M_p,Z)\cong\pi_{2n-1}Z$. Rationalizing \eqref{common:primary-obstruction} and identifying the homology via
projection onto the first factor gives the obstruction homomorphism $k$ and the formula $\ker k=\im g_*$.

If $f$ lifts, then $f_*$ factors through $g_*$. Conversely, assume the image condition. Adjoin
constant higher coordinates to $f$ using \eqref{even:postnikovproduct}. The obstruction evaluates to
$k\circ f_*=0$. Since $\pi_{2n-1}Z$ is divisible, the $\operatorname{Ext}$ term vanishes. As $\dim
X\le2n$, we obtain a lift to $Z$. The homotopy fiber of $g$ is $n$-connected and $\dim X_0\le n$, so
the restriction of this lift is homotopic to the prescribed lift over $f|_{X_0}$. Relative homotopy
lifting gives a lift agreeing with the prescribed one on $X_0$.
\end{proof}

\section{The main theorems in odd and even dimensions}\label{sec:odd}\label{sec:even}

In odd dimensions a family of sections parametrized by a rational sphere gives a coherent family of
integral periods. In even dimensions we use the square-zero relation for the transverse volume class
to determine when periods vanish. To prove surjectivity we construct sections from cohomology classes
of intermediate degree. This construction also realizes individual real periods in odd dimensions.

\subsection{Coherent families of integral periods in odd dimensions}\label{subsec:odd-periods}

Let $n\ge3$ be odd. We use CW replacements or singular Kan models for derived mapping spaces as needed.

\begin{lemma}[Maps from a rational sphere]\label{odd:tel}\label{odd:map-rational}
For every $W\in\R$ there is a based map $\theta_W:K_{\Q,n}\to F_n$ such that evaluation of $\theta_W^*u_F$ on $H_n(K_{\Q,n};\Z)=\Q$ is $t\mapsto Wt$. Moreover, for any based
space $F$ each path component of $\Map_*(K_{\Q,n},F)$ is weakly equivalent to a simple rational space. Here a space is
simple if its fundamental group is abelian and acts trivially on its higher homotopy groups.
\end{lemma}

\begin{proof}
Since $n$ is odd, \eqref{common:em-homology} and the fact that the positive-degree integral homology groups
are rational vector spaces identify $K_{\Q,n}$ as a model for the Moore space $M(\Q,n)$.
By the universal coefficient sequence for homotopy with coefficients \cite[Proposition~4H.2]{Hatcher},
every homomorphism $\Q\to\pi_nF_n$ is realized by a based map $K_{\Q,n}\to F_n$.
Choose the homomorphism corresponding to $t\mapsto Wt$ under \eqref{common:localhomotopy}.
Naturality of the Hurewicz isomorphism gives the required evaluation.

Uniqueness of Moore spaces \cite[Example~4.34]{Hatcher} gives $K_{\Q,n}\simeq\Sigma M(\Q,n-1)$. Thus $\Map_*(K_{\Q,n},F)$ is weakly equivalent to a loop space, and loop
multiplication identifies each component up to homotopy with the component of the constant map. This
component is a connected $H$-space and hence simple. On $\pi_i=[S^i\wedge K_{\Q,n},F]_*$ for $i\ge1$, multiplication by $m$ is
precomposition with the smash product of the degree-$m$ map on $S^i$ and the identity of $K_{\Q,n}$. The only
nonzero reduced integral homology of the simply connected space $S^i\wedge K_{\Q,n}$ is $\widetilde H_{n+i}=\Q$. This self-map acts there
by multiplication by $m$ and is a homotopy equivalence for $m\ge1$ \cite[Corollary~4.33]{Hatcher}. Hence each $\pi_i$ is uniquely
divisible, and a CW model of each component is rational \cite[Theorem~6.1.1]{MayPonto}.
\end{proof}

To choose maps from the rational sphere fiber by fiber, we use a point-evaluation fibration. Let $F\to E\xrightarrow{\nu}B$
be a fibration of simply connected spaces and $\theta:K_{\Q,n}\to F$ a based map. Define the space of maps into
individual fibers and its evaluation by
\begin{equation}\label{odd:eq:evaluation}
 \mathcal E_{K_{\Q,n}}\coloneqq\{(b,f)\mid b\in B,\ f:K_{\Q,n}\to\nu^{-1}(b)\},\qquad
 \epsilon(b,f)\coloneqq f(0).
\end{equation}
Give $\mathcal E_{K_{\Q,n}}$ the compactly generated subspace topology from $B\times\Map(K_{\Q,n},E)$.
Take $K_{\Q,n}$ to be the geometric realization of a simplicial abelian group, with $0$ a $0$-cell.
Given a homotopy of the evaluation points in $E$, project it to $B$ by $\nu$ and prescribe the original
homotopy at $0\in K_{\Q,n}$. The relative homotopy lifting property of $\nu$ then shows that
$\epsilon:\mathcal E_{K_{\Q,n}}\to E$ is a fibration. Since $E$ is simply connected, its monodromy action on the components of the fiber is
trivial. Let $b_0$ be the basepoint of $B$ and $i_F:F=\nu^{-1}(b_0)\hookrightarrow E$ the inclusion. Restricting to the path component $\mathcal E_{K_{\Q,n},\theta}$
containing $(b_0,i_F\theta)$ gives
\begin{equation}\label{odd:eq:evaluation-component}
 N_\theta\coloneqq\Map_*(K_{\Q,n},F)_\theta\longrightarrow\mathcal E_{K_{\Q,n},\theta}
 \xrightarrow{\epsilon}E.
\end{equation}
Here $N_\theta$ is the fiber component containing $\theta$.

\begin{proposition}[Rationalization of the evaluation fibration]\label{odd:descent}
In the chosen homotopy models there is a weak equivalence over $E$
\begin{equation}\label{odd:eq:cartesian}
 \mathcal E_{K_{\Q,n},\theta}\simeq_E E\times^h_{E_\Q}(\mathcal E_{K_{\Q,n},\theta})_\Q
\end{equation}
where the homotopy pullback is taken along the rationalization $E\to E_\Q$. Let $\sigma,h:X\to E$ be maps from a
finite CW complex $X$ such that $\ell\sigma\simeq\ell h$. Then $\sigma^*\mathcal E_{K_{\Q,n},\theta}$ admits a section
if and only if $h^*\mathcal E_{K_{\Q,n},\theta}$ does.
\end{proposition}

\begin{proof}
By Lemma~\ref{odd:map-rational}, $N_\theta$ is simple and rational. Since the base $E$ is simply connected, the total space of \eqref{odd:eq:evaluation-component}
is nilpotent \cite[Lemma~4.7]{BerglundMadsen}. The rationalization theorem for fibrations with simply connected base \cite[Theorem~6.2.6 and the remarks immediately following it]{MayPonto} identifies
the map from $N_\theta$ to the homotopy fiber of the rationalized projection as a rationalization. Since $N_\theta$
is already rational, this map is a weak equivalence, giving \eqref{odd:eq:cartesian}.

Suppose $h$ has a lift $\widetilde h$. The map $\sigma$, the composite $\ell\circ\widetilde h$ with the rationalization of the total space, and the
homotopy $\ell\sigma\simeq\ell h$ define a map from $X$ to the homotopy pullback in \eqref{odd:eq:cartesian}.
Since $X$ is a CW complex, the weak equivalence gives a lift of $\sigma$ up to homotopy.
The homotopy lifting property of $\epsilon$ gives a strict lift. Interchanging $\sigma$ and $h$ proves the converse.
\end{proof}

\begin{proposition}[A family of sections parametrized by a rational sphere]\label{odd:section-map}
Under the hypotheses of Theorem~\ref{thm:odd}, there is a based map $A_M:K_{\Q,n}\to\mathcal S_M$
into the reference component of the section space of $\tau_M^*E_n\to M$ such that
\begin{equation}\label{odd:eq:section-period}
 \langle A_M^*\ev_{a_0}^*u_F,t\rangle=\totalVol\,t\qquad(t\in\Q).
\end{equation}
\end{proposition}

\begin{proof}
Choose a based map $\kappa_M:(M,a_0)\to(K_{\Q,n},0)$ classifying the rational orientation class that pairs
to one with $[M]$. Thus $(\kappa_M)_*[M]=1$. Take $\theta=\theta_{\totalVol}$ from Lemma~\ref{odd:tel} and set
$h\coloneqq i_F\theta\kappa_M:M\to E_n$. Both $\sigma_0$ and $h$ pull the volume class back to $[v]$.
The tangent coordinate of $h$ is constant, and $\ell\tau_M$ is nullhomotopic by \eqref{odd:eq:bso} and rational Pontryagin-triviality.
$\ell\sigma_0$ and $\ell h$ have homotopic composites with the map in \eqref{common:joint}. The fiber of $(\nu_\Q,u_\Q)$ is the homotopy fiber of
the volume classifying map $(F_n)_\Q\to K(\R,n)$. By \eqref{common:localhomotopy} this map induces an isomorphism on $\pi_n$, and its target has
zero $\pi_{n+1}$, so its fiber is $n$-connected. We may therefore lift the homotopy in the target relative to
its two ends over the $n+1$-dimensional pair $(M\times[0,1],M\times\{0,1\})$.

The assignment $x\mapsto(z\mapsto i_F\theta(z+\kappa_M(x)))$ defines a section of $h^*\mathcal E_{K_{\Q,n},\theta}$. Its value at $z=0$ is $h(x)$, and at $x=a_0$ it gives $i_F\theta$, so the image
lies in the specified path component. Apply Proposition~\ref{odd:descent} to this section and the homotopy $\ell\sigma_0\simeq\ell h$ just
obtained. The resulting section of $\sigma_0^*\mathcal E_{K_{\Q,n},\theta}$ corresponds to a family $f_x:K_{\Q,n}\to\nu^{-1}(\tau_M(x))$ satisfying $f_x(0)=\sigma_0(x)$. Its adjoint $A_M(z)(x)=f_x(z)$
thus satisfies $A_M(0)=s_0$. Since $K_{\Q,n}$ is path connected, the image lies in $\mathcal S_M$. The points corresponding to
$f_{a_0}$ and $i_F\theta$ belong to the same path component of the evaluation space, so the maps are homotopic as
maps to $E_n$. In particular they have the same pullback of $u$. Under the volume frame at $a_0$, the map
$\ev_{a_0}A_M$ corresponds to $f_{a_0}$. This proves \eqref{odd:eq:section-period}.
\end{proof}

\begin{proof}[Proof of Theorem~\ref{thm:odd}]
Since $K_{\Q,n}$ is simply connected, $A_M$ lifts to a based map $\widetilde A_M$ into $\widetilde{\mathcal S}_M$. With integral coefficients put
\begin{equation}\label{odd:eq:alpha}
 \alpha_{\cH}^{\Q}\coloneqq \mathsf T_{\cH}(\widetilde A_M)_*:
 \Q=H_n(K_{\Q,n};\Z)\longrightarrow H_n(B\cH^\delta;\Z)
\end{equation}
and define $\alpha_{\cG}^{\Q}=i_*\alpha_{\cH}^{\Q}$. The $\cH$ formula in Proposition~\ref{common:comparison} and the evaluation \eqref{odd:eq:section-period} give $P_{\cH}(\alpha_{\cH}^{\Q}(t))=\totalVol\,t$. Naturality
under inclusion gives $P_{\cG}(\alpha_{\cG}^{\Q}(t))=\totalVol\,t$ as well. Since $\totalVol>0$, both maps are injective.

Corollary~\ref{common:lone} gives
\begin{equation}\label{odd:eq:norm}
 \|\alpha_{\cG}^{\Q}(t)_\R\|_1\ge\|M\|\,|t|
\end{equation}
and the same lower bound holds on $\cH$. The images are divisible, so they are injective objects in the category
of abelian groups. Each inclusion therefore splits.
\end{proof}

\subsection{The even-dimensional period criterion and intermediate cohomology}\label{subsec:even-periods}

First suppose $n=2m\ge4$.

\begin{proposition}[A quadratic identity for the top volume class]\label{even:quadratic}
For each group $\cL$ of Section~\ref{subsec:flat} we have
\begin{equation}\label{even:quadraticformula}
 2\,\totalVol\,\xi_{\cL}+(\pi_{\cL})_!(\Theta_{\cL,\mathrm{mix}}^2)=0.
\end{equation}
For $\cL=\cG$, the mixed components of fiber degrees $1,n-1$ make no contribution to the second term after
restriction to $B\cH^\delta$.
\end{proposition}

\begin{proof}
Expand $(\pi_{\cL})_!(\Theta_{\cL}^2)=0$ using \eqref{common:kun}. Since $n$ is even, the two products of the pure base and pure fiber components
give $2\,\totalVol\,\xi_{\cL}$. Of the remaining products, only those between mixed components whose fiber degrees sum to $n$
contribute.

For $\cL=\cG$, restrict the flat bundle to the loop in $B\cG^\delta$ represented by $g\in\cG$.
Choose a volume-preserving isotopy $g_t$ from the identity to $g$, constant near both endpoints, and let $X_t$
be its generating vector field. The map $(t,x)\mapsto[t,g_t x]$ trivializes the bundle with gluing
$(1,x)\sim(0,g^{-1}x)$. In this trivialization, uniqueness of the projection gives $r\simeq\pr_M$.
With the base oriented first, integration of the transverse volume form in the interval direction gives
$[\int_0^1g_t^*(\iota_{X_t}v)\,dt]=\Flux(g)$ \cite[Lemma~8, Proof of Theorem~10]{KKM}.
Thus the component of bidegree $(1,n-1)$ vanishes on $\cH$, as does its product with the component of fiber degree $1$.
\end{proof}

If $b_{[1,m]}(M)=0$, Poincar\'e duality gives $H^j(M;\R)=0$ for every $0<j<n$, so \eqref{even:quadraticformula}
gives $\xi_{M,n}=0$. If $b_{[2,m]}(M)=0$, cohomology vanishes in degrees $2,\ldots,n-2$.
The only possible contributions come from the paired components of fiber degrees $1$ and $n-1$,
whose product vanishes on $\cH$. Thus $\xi_{M,n}|_{\cH}=0$ in this case.

To prove surjectivity we construct families of sections parametrized by products of spheres.
We begin with a cohomology class that we will use in both odd and even dimensions.

\begin{lemma}\label{even:two-direction}
Let $M$ be a closed connected oriented $n$-manifold with $n\ge3$, and let $\mu_M$ be the rational orientation
class satisfying $\langle\mu_M,[M]\rangle=1$. Suppose $a\in H^p(M;\Q)$ and $b\in H^{n-p}(M;\Q)$ satisfy $ab=\mu_M$,
where $1\le p\le n-1$. On $Y=S^{n-p}\times S^p$, let $x,y$ be the normalized orientation classes pulled back
from the first and second spheres, respectively. Give $X=Y\times M$ the product orientation and set
\[
 A\coloneqq 1\times\mu_M+x\times a,\qquad
 B\coloneqq xy\times1+(-1)^{n+p+1}y\times b.
\]
For all $V,W\in\R$, the map $X\to K(\R,n)$ classifying $VA+WB$ induces the zero map on $H_{2n}(-;\Q)$.
The restrictions of $VA+WB$ to $\{*\}\times M$ and to $Y\times\{a_0\}$, for $a_0\in M$, are $V\mu_M$ and $Wxy$, respectively.
\end{lemma}

\begin{proof}
Degrees and the identities $x^2=y^2=0$ give $A^2=B^2=0$. The sign rule for cross products and $ab=\mu_M$ give $AB=(-1)^nxy\times(\mu_M-ab)=0$. Let $(A,B):X\to K(\Q^2,n)$
classify the two rational classes. By the rational cohomology calculation for $K(\Q,n)$ \cite[Proposition~6.7.2]{MayPonto} and the K\"unneth
theorem, $H^{2n}(K(\Q^2,n);\Q)$ is spanned by the product of the two standard degree-$n$ classes and, when $n$ is even,
their squares. Every generator pulls back to zero. Naturality of the rational evaluation pairing
therefore shows that $(A,B)$ induces zero on $H_{2n}(-;\Q)$. The map classifying $VA+WB$ is its composite with $K(\Q^2,n)\to K(\R,n)$
induced by $(r,s)\mapsto Vr+Ws$. The two restrictions follow from the definitions of $A,B$.
\end{proof}

In even dimension $n=2m\ge4$, form the homotopy pullback along the inclusion $j_e$ of the Euler factor in \eqref{even:bso},
\begin{equation}\label{even:eulerpullback}
 E'_n\coloneqq K_{\Q,n}\times_{B\SO(n)_\Q}(E_n)_\Q
\end{equation}
using the fibration model of $\nu_\Q$. Write $e,u$ for the Euler and transverse volume classes pulled back to
$E'_n$. By the fiber sequence $(F_n)_\Q\to E'_n\to K_{\Q,n}$ and \eqref{common:localhomotopy}, the space $E'_n$ is $(n-1)$-connected and rational, and
\begin{equation}\label{even:lowpi}
 (e,u)_*:\pi_n E'_n\xrightarrow{\ \cong\ }\Q\oplus\R.
\end{equation}

\begin{proposition}[Families of sections from intermediate cohomology]\label{even:family}
Suppose that $M$ satisfies the hypotheses of Section~\ref{sec:intro}, that $n\ge3$, and that $M$ is rationally Pontryagin-trivial.
If $H^p(M;\Q)\ne0$ for some $1\le p\le n-1$, then for each $W\in\R$ there exists a based map
$a_W:Y\to\mathcal T_M$, with $Y=S^{n-p}\times S^p$, such that
\begin{equation}\label{even:familyperiod}
 \left\langle a_W^*\ev_{a_0}^*(u_F)_\Q,[Y]\right\rangle=W.
\end{equation}
If $2\le p\le n-2$, then $a_W$ lifts to $\widetilde{\mathcal T}_M$.
\end{proposition}

\begin{proof}
Choose $a,b$ by Poincar\'e duality as in Lemma~\ref{even:two-direction}. With $A,B$ as in that lemma, set
$U_W\coloneqq\totalVol\,A+WB$. Let $X=Y\times M$, $X_0=\{*\}\times M$, and let $\pr_M:X\to M$ be the projection.
We use the same notation for cohomology classes and their classifying maps.
By rational Pontryagin-triviality and \eqref{odd:eq:bso}--\eqref{even:bso}, the rationalized classifying map
of $TM$ is nullhomotopic in odd dimensions and homotopic to $j_e\circ e(TM)$ in even dimensions.
Transport the reference section along this homotopy, keeping the volume coordinate fixed.
In each case we seek a lift of the map $f$ along the map $g$ indicated below.
\[
\begin{array}{c@{\qquad}l@{\qquad}l}
 n & g & f\ \text{(domain $X$)} \\
 \text{odd} & (u_F)_\Q:(F_n)_\Q\to K(\R,n) & U_W \\
 \text{even} & (e,u):E'_n\to K_{\Q,n}\times K(\R,n) & (\pr_M^*e(TM),U_W)
\end{array}
\]
By \eqref{common:localhomotopy} and \eqref{even:lowpi}, in both cases $g$ is the projection to the $n$-th Postnikov stage. Since $U_W|_{X_0}=[v]$, the volume
classifying maps are homotopic on $X_0$. Use homotopy extension for the CW pair $(X,X_0)$ to adjust the volume
coordinate of $f$ so that the transported reference section is a prescribed lift of $f|_{X_0}$. Since $\dim X=2n$,
$\dim X_0=n$, and $H_{2n}(X;\Q)=\Q\cdot[X]$, Lemma~\ref{even:postnikov} reduces the problem to showing $f_*[X]\in\im g_*$.

By Lemma~\ref{even:two-direction}, the pure volume component of $f_*[X]$ vanishes. In odd dimensions this gives
$f_*[X]=0$, so the image condition holds. In even dimension $n=2m$, the Euler coordinate factors through
$M$, and the pure Euler component also vanishes. The mixed component lies in $\Q\otimes_\Q\R\cong\R$.
Since $e(TM)=\chi(M)\mu_M$ and $U_W|_{Y\times\{a_0\}}=Wxy$, it is given by
\[
 t\coloneqq\langle\pr_M^*e(TM)\smile U_W,[X]\rangle=\chi(M)W.
\]
In the order determined by the pure Euler, mixed, and pure volume components, $f_*[X]=(0,t,0)$.
It remains to construct a class mapping to this triple.

On $E_2$, let $F$ be the homotopy fiber of the map $(e_2,u_2):E_2\to K(\Z,2)\times K(\R,2)$ classifying the Euler class $e_2$ of the normal bundle
and the transverse area class $u_2$. By Nariman \cite[preprint v2, proof of Theorem~3.9, pp.~28--29]{Nariman}, $F$ is $2$-connected. Under a suitable normalization
of $\pi_3F\cong\R$, the class in $H^3(F;\R)$ corresponding to the identity homomorphism transgresses to the square of the
universal class of the second factor. Up to sign, this transgression is the primary obstruction in
degree four. On $S^2\times S^2$ write $\alpha,\beta$ for the pullbacks of the positive integral generators from the first
and second spheres, respectively. For any $s\in\R$ the obstruction pulled back along $(\alpha,s\beta)$ is $(s\beta)^2=0$. The
domain is four-dimensional, so there is a lift $S^2\times S^2\to E_2$.

Take the external product of these lifts with $s=t$ on the first pair and $s=1$ on the remaining $m-1$ pairs.
Products of area-preserving germs preserve volume, giving a map from $P=(S^2\times S^2)^m$ to $E_n$. Let $\alpha_j,\beta_j$ be the
pullbacks from the $j$-th pair. The normal bundle is the external direct sum of the oriented plane
bundles, with total Pontryagin class $\prod_j(1+\alpha_j^2)=1$. Choosing nullhomotopies for the Pontryagin coordinates in \eqref{even:bso}
gives, by the homotopy pullback \eqref{even:eulerpullback}, a lift $k:P\to E'_n$ of $P\to E_n\to(E_n)_\Q$ up to homotopy. The product formula for Euler
classes and the exterior product of the transverse area forms give
\[
 k^*e=\prod_{j=1}^m\alpha_j,\qquad k^*u=t\prod_{j=1}^m\beta_j .
\]
Up to homotopy, the two coordinate maps thus factor through the product of the first spheres and the
product of the second spheres, respectively. Both products have dimension $n$, so both pure components
vanish. With the product orientation, $\langle k^*(eu),[P]\rangle=t$. Thus $g_*k_*[P]=(0,t,0)=f_*[X]$.

Lemma~\ref{even:postnikov} now extends the prescribed lift on $X_0$ over $X$. Returning to the original tangent coordinate
by Lemma~\ref{common:tangent-transport} gives a based family of sections $a_W:Y\to\mathcal T_M$. The corresponding lift on $Y\times M$ has volume class $U_W$ for
the chosen dual classes $a,b$. The pullback of the point-evaluation class is $U_W|_{Y\times\{a_0\}}=Wxy$, which evaluates on $[Y]$
to $W$. If $2\le p\le n-2$, then $Y$ is simply connected, and $a_W$ lifts to the universal cover $\widetilde{\mathcal T}_M$.
\end{proof}

We next pass from rational homology classes on $\mathcal T_M$ or its universal cover to integral periods on
the corresponding discrete group. The argument applies in every dimension $n\ge3$.

A Polish space is a separable, completely metrizable space, and a Polish group is a topological group
whose underlying space is Polish. A subset of a Polish space is analytic if it is the Borel image of
a Borel subset of a Polish space.

\begin{lemma}\label{even:analytic}
Let $L$ be a Polish group and $c:L^d\to\R$ a Borel group cocycle of degree $d\ge1$. Its integral period group $\operatorname{Per}(c)\coloneqq\{\langle[c],z\rangle\mid z\in H_d(BL^\delta;\Z)\}$ is
an analytic subgroup of $\R$. If its rational linear span is $\R$, then $\operatorname{Per}(c)=\R$.
\end{lemma}

\begin{proof}
In the unnormalized bar complex write integral chains as finite sums with coefficients $\pm1$, allowing
repetitions. For a fixed number $N$ of terms and fixed signs, these chains are parametrized by $L^{dN}$. The
boundary vanishes precisely when its positive and negative terms can be paired with equal labels. The
labels vary continuously and there are only finitely many pairings, so the parameter set of cycles is
closed. Evaluation by $c$ is Borel, and the countable union of its images over all lengths and sign
choices is $\operatorname{Per}(c)$. Addition and negation of chains make it a subgroup. Being analytic, it has the Baire
property \cite[Theorem~27]{Moore}.

If the rational span is $\R$, then $\R=\bigcup_{N\ge1}N^{-1}\operatorname{Per}(c)$, so the Baire category theorem shows that $\operatorname{Per}(c)$ is nonmeager. By
Pettis's lemma and the subgroup property, $\operatorname{Per}(c)$ contains a neighborhood of zero \cite[Lecture~13, Lemma~6]{Moore}. It is therefore an
open subgroup of $\R$, and $\operatorname{Per}(c)=\R$.
\end{proof}

\begin{lemma}\label{even:periodtransfer}
Suppose the point-evaluation class on $\mathcal T_M$ takes every real value on $H_n(\mathcal T_M;\Q)$.
Then $\im P_{\cG}=\R$. If its pullback to $\widetilde{\mathcal T}_M$ takes every real value on
$H_n(\widetilde{\mathcal T}_M;\Q)$, then $\im P_{\cH}=\R$.
\end{lemma}

\begin{proof}
Fix $r\in\R$. In the first case choose a rational homology class of period $r$ on $\mathcal T_M$,
and in the second case choose one on $\widetilde{\mathcal T}_M$.
The rational homology isomorphism of Section~\ref{subsec:section-models} gives a preimage on $\mathcal S_M$
or on its universal cover, respectively. This preimage still has period $r$. Since
$H_n(-;\Q)=H_n(-;\Z)\otimes\Q$, some positive integer multiple of the preimage comes from integral homology.
Applying $\mathsf T_{\cL}$ and Proposition~\ref{common:comparison}, with $\cL=\cG$ in the first case and $\cL=\cH$ in the second,
gives $Nr\in\im P_{\cL}$ for some integer $N>0$. Thus $\im P_{\cL}$ spans $\R$ over $\Q$.

Since $M$ is compact, $\Diff(M)$ is open in the Polish space $C^\infty(M,M)$.
The subgroup preserving $v$ is closed in $\Diff(M)$, and its identity component $\cG$ is closed in that subgroup.
It follows that $\cG$ is Polish. The continuity of flux implies that $\cH$ is Polish as well.
Let $c$ be the inhomogeneous cocycle on $\cL$ obtained from the continuous representative of Proposition~\ref{common:continuous}.
Then $\operatorname{Per}(c)=\im P_{\cL}$, so Lemma~\ref{even:analytic} gives $\im P_{\cL}=\R$.
\end{proof}

\begin{proposition}[Intermediate cohomology and realization of all real periods]\label{oddreal:main}
Suppose that $M$ satisfies the hypotheses of Section~\ref{sec:intro}, that $n\ge3$, and that $M$ is rationally Pontryagin-trivial.
If $H^p(M;\Q)\ne0$ for some $1\le p\le n-1$, then $\im P_{\cG}=\R$.
If one can choose $p$ with $2\le p\le n-2$, then $\im P_{\cH}=\R$ as well.
\end{proposition}
\begin{proof}
By Proposition~\ref{even:family} every $W\in\R$ occurs as a point-evaluation period on $\mathcal T_M$, and for $2\le p\le n-2$ it occurs on the
universal cover as well. Lemma~\ref{even:periodtransfer} gives the required integral periods.
\end{proof}

\begin{proof}[Proof of Theorem~\ref{thm:even}]
Vanishing follows from Proposition~\ref{even:quadratic} and the argument immediately following its proof. For
surjectivity suppose $M$ is rationally Pontryagin-trivial. If $b_{[1,m]}(M)>0$, there is a degree $p$ in the range $1\le p\le m$
with $H^p(M;\Q)\ne0$. If $b_{[2,m]}(M)>0$, such a degree lies in the range $2\le p\le m\le n-2$. Proposition~\ref{oddreal:main} gives $\im P_{\cG}=\R$ and $\im P_{\cH}=\R$, respectively.
\end{proof}

\subsection{The criterion under nonpositive curvature}

\begin{proposition}[Vanishing on nonpositively curved manifolds]\label{scope:bott}
Let $M^n$ be a closed connected oriented smooth manifold with centerless fundamental group.
Suppose that $M$ admits a metric of nonpositive sectional curvature and that $p_i(TM)\ne0$ in $H^{4i}(M;\Q)$
for some $i>0$. Then for every positive volume form and every subgroup $\cL$ as in Section~\ref{subsec:flat}, we have $\xi_{\cL}=0$.
\end{proposition}

\begin{proof}
Take any $z\in H_n(B\cL^\delta;\Z)$. Since $z$ comes from a finite subcomplex, apply Thom's realization theorem \cite[Théorème~III.4]{ThomRealization} to a finite
polyhedral model for that subcomplex. This yields a closed oriented smooth $n$-manifold $N$ and a map
$f:N\to B\cL^\delta$ with $f_*[N]=kz$ for some positive integer $k$. Realize the pulled-back flat bundle $\pi_f:E_f\to N$ smoothly as the
suspension of its monodromy, and write $\Theta_f$ for its transverse volume class. Identify the normal
bundle of the horizontal foliation with the vertical tangent bundle $T_vE_f$. In flat bundle coordinates
$U\times M$, keeping vectors of $T_yM$ constant along the plaques $U\times\{y\}$ defines leafwise parallel transport. A
connection on the normal bundle extending this partial connection is a Bott connection. Applying
Hurder \cite[\S1, equation~(1.2)]{Hurder} to $c_{2i}$ of the complexified normal bundle, we see that for $i>0$ every term in the
Chern--Weil form representing $p_i(T_vE_f)=(-1)^i c_{2i}(T_vE_f\otimes\mathbb C)$ contains at least one transverse differential $dy^j$. The transverse
volume form representing $\Theta_f$ is locally $g(y)\,dy^1\wedge\cdots\wedge dy^n$. The wedge product of the two forms therefore vanishes,
giving $\Theta_f\smile p_i(T_vE_f)=0$.

The projection of Section~\ref{subsec:flat} gives a fiber homotopy equivalence $(\pi_f,r_f):E_f\to N\times M$.
The nonvanishing of $p_i(TM)$ implies $n\ge4i\ge4$. Since $N$ admits a finite triangulation,
Bustamante, Farrell, and Jiang \cite[Theorem~B]{BFJ} identify $T_vE_f$ and $r_f^*TM$ as topological $\R^n$-bundles.
Their rational Pontryagin classes therefore agree \cite[Proof of Corollary~B.1]{BFJ}.
By Poincar\'e duality, choose $\beta\in H^{n-4i}(M;\R)$ with $p_i(TM)\smile\beta=[v]$. Then
\[
 \Theta_f\smile r_f^*[v]
  =\Theta_f\smile p_i(T_vE_f)\smile r_f^*\beta=0.
\]
The fiber integration formula of Section~\ref{subsec:flat} gives $f^*\xi_{\cL}=0$. Hence $k\langle\xi_{\cL},z\rangle=\langle f^*\xi_{\cL},[N]\rangle=0$. Since $z$ was arbitrary, the universal
coefficient theorem over the reals gives $\xi_{\cL}=0$.
\end{proof}

Combining Proposition~\ref{scope:bott} with the main theorems gives the following criterion whenever $M$ admits
a metric of nonpositive sectional curvature. In even dimension $n=2m\ge4$, take $(\cL,a)$ to be either $(\cG,1)$ or $(\cH,2)$.
Then $\im P_{\cL}=\R$ if and only if $M$ is rationally Pontryagin-trivial and $b_{[a,m]}(M)>0$.
In odd dimension $n\ge3$, for each of $\cG$ and $\cH$ there is a nonzero integral period if and only if
$M$ is rationally Pontryagin-trivial (Theorem~\ref{thm:odd}).

\subsection{The bounded volume class of a closed real hyperbolic manifold}\label{subsec:bm-proof}

A closed connected oriented real hyperbolic manifold $M$ of dimension $n\ge3$ is aspherical and has
centerless fundamental group. Hyperbolic holonomy identifies $\Gamma$ with a torsion-free cocompact lattice
in $\operatorname{SO}_0(n,1)$. It follows from Tshishiku's theorem \cite[Theorem~6.7]{Tshishiku} that
$p_i(TM)=0\in H^{4i}(M;\Q)$ for all $i>0$. Thus Theorems~\ref{thm:odd} and~\ref{thm:even} apply to any
positive volume form in odd and even dimensions, respectively. When $n=2m$, the Gauss--Bonnet formula gives
$\operatorname{sgn}\chi(M)=(-1)^m$ \cite[equation~(0.1)]{KellerhalsZehrt}.

\begin{proof}[Proof of Theorem~\ref{thm:bm-detection}]
Use the homomorphism $J:\cG\to H$ of Proposition~\ref{common:bm-comparison} and let $v$ be the hyperbolic volume form. In odd dimensions
put $y_t=J_*\alpha_{\cG}^{\Q}(t)$. Theorem~\ref{thm:odd} gives $P_{\cG}(\alpha_{\cG}^{\Q}(t))=\totalVol\,t$, so \eqref{common:bm} yields $\langle \operatorname{comp}(\mathcal B_M),(y_t)_\R\rangle=-\totalVol^2t$. Thus the additive map $t\mapsto y_t$ is injective. Its image is
divisible and hence a direct summand.

In even dimension $n=2m$, Poincar\'e duality shows that $b_{[1,m]}(M)=0$ if and only if $M$ is a rational homology
sphere. An even-dimensional rational homology sphere has $\chi(M)=2$. When $m$ is odd, $\chi(M)<0$ excludes this
possibility. When $m$ is even, it is excluded by hypothesis. Thus $P_{\cG}$ is surjective by Theorem~\ref{thm:even}. By
\eqref{common:bm} the value on $J_*z$ is $\totalVol\,P_{\cG}(z)$, so every real number occurs.

In either case we can choose an integral class $z$ with $P_{\cG}(z)=\totalVol$. The ordinary class
$\operatorname{comp}(\mathcal B_M)$ pairs nontrivially with $J_*z$. The duality inequality gives
$0<\totalVol^2\le\|J^*\mathcal B_M\|_\infty\,\|z_\R\|_1$. Since pullback does not increase the seminorm, $\|\mathcal B_M\|_\infty>0$.
\end{proof}

For an even-dimensional rational homology sphere, Theorem~\ref{thm:even} and \eqref{common:bm} give $J^*\operatorname{comp}(\mathcal B_M)=0$.

\section{Coherent families of periods in dimension three}\label{sec:threeD-refinement}

Throughout this section the three-manifold $M$ satisfies the standing assumptions. We characterize the
normalized integral periods on the flux kernel using the sixth integral homology of $F_3$, the
classifying space for framed Haefliger structures with transverse volume. We use the Sol group to
construct a coherent family over $E_{\log}$.

\subsection{Sixth homology and coherent families of periods}
We first determine when a map $S^3\vee S^3\to F_3$ extends over $S^3\times S^3$.
We then use $H_6(F_3;\Z)$ to construct additive lifts of normalized periods to integral homology
(Theorem~\ref{h3:geo:integer-additive-section}).

Write $F=F_3$. Since oriented three-manifolds are parallelizable, we can choose a tangent framing
compatible with $v$ that agrees with the chosen frame at $a_0$. In this framing the reference component
$\mathcal S_M$ of McDuff's section space has the mapping-space model $\Map(M,F)_\sigma$, where $\sigma:M\to F$ corresponds to the reference
section and satisfies $\sigma^*u_F=[v]$.

Put $K_{\R,3}\coloneqq K(\R,3)$, and let $\iota\in H^3(K_{\R,3};\R)$ be the universal class corresponding to the identity homomorphism of the
coefficient group $\R$. Let $q:F\to K_{\R,3}$ be the based classifying map of the volume class, so that $q^*\iota=u_F$. By \eqref{common:localhomotopy}, the
space $F$ is 2-connected, $q$ induces an isomorphism on $\pi_3$, and $\pi_4F=0$. Moreover,
\[
 H_6(K_{\R,3};\Z)\cong\R\wedge_{\Q}\R,\qquad H_7(K_{\R,3};\Z)=0.
\]
Here $\R\wedge_{\Q}\R$ denotes the exterior square of $\R$ over $\Q$. These identifications follow from
\eqref{common:em-homology}, since the integral homology groups in positive degrees are rational vector spaces.
We normalize the first isomorphism so that the Pontryagin product of two spherical classes in degree three,
of volumes $s$ and $t$ in that order, corresponds to $s\wedge t$.

Write $\smile_\wedge$ for the cup product followed by the coefficient map $\R\otimes_\Z\R\to\R\wedge_{\Q}\R$, $s\otimes t\mapsto s\wedge t$. For the map $\varphi:X\to K_{\R,3}$ classifying $u\in H^3(X;\R)$
and for $z\in H_6(X;\Z)$, we have
\begin{equation}\label{h3:src:finite-wedge}
 \varphi_*z=\tfrac12\langle u\smile_\wedge u,z\rangle .
\end{equation}
Indeed, on the Pontryagin product of spherical classes of volumes $s$ and $t$, the two cross terms in
$\iota\smile_\wedge\iota$ give $s\wedge t$ and $-t\wedge s$. Their sum is $2(s\wedge t)$.
Since these Pontryagin products generate $H_6(K_{\R,3};\Z)$, the formula follows by naturality.

Define
\[
 I=\im\bigl(q_*:H_6(F;\Z)\to\R\wedge_{\Q}\R\bigr),\qquad
 I_{\mathrm{div}}=\{q_*\gamma(1):\gamma\in\Hom_\Z(\Q,H_6(F;\Z))\}
\]
and define the two candidate sets of normalized periods by
\[
 \mathcal P=\{u\in\R:1\wedge u\in I\},\qquad
 \mathcal P_{\mathrm{div}}=\{u\in\R:1\wedge u\in I_{\mathrm{div}}\}.
\]
Thus $w\in I_{\mathrm{div}}$ precisely when there is a homomorphism $\gamma:\Q\to H_6(F;\Z)$
with $q_*\gamma(r)=rw$ for every $r\in\Q$.

For $\lambda>0$, the transverse dilation $h_\lambda$ of Section~\ref{subsec:section-models} acts on $F=F_3$.
The corresponding map $T_\lambda(s\wedge t)=\lambda s\wedge\lambda t$ on $H_6(K_{\R,3};\Z)$ preserves both $I$ and $I_{\mathrm{div}}$.
For each positive integer $d$, the operator $dT_{1/d}$ is multiplication by $1/d$, so $I$ is a $\Q$-vector subspace.
For $I_{\mathrm{div}}$, closure under rational scalars also follows by precomposing $\gamma$ with the corresponding
scalar multiplication on $\Q$.

Denote the homotopy fiber of $q$ by $Z_q$.

\begin{lemma}\label{h3:geo:primary-obstruction}
Let $s,t:S^3\to F$ be based maps of volumes $W=\langle s^*u_F,[S^3]\rangle$ and $V=\langle t^*u_F,[S^3]\rangle$. Then $s\vee t:S^3\vee S^3\to F$ extends over $S^3\times S^3$ if and only if $W\wedge V\in I$.
\end{lemma}
\begin{proof}
In the topological abelian group model of $K_{\R,3}$ put $f(x,y)=q(s(x))+q(t(y))$. On the wedge, $s\vee t$ is a prescribed lift of $f$.
With the product oriented in the indicated order, $f_*[S^3\times S^3]=W\wedge V$.
Write $M_q$ for the mapping cylinder of $q$. The fiber $Z_q$ is $4$-connected, and $S^3\times S^3$ has a single
cell relative to the wedge, in degree six. Under the relative Hurewicz isomorphism of Section~\ref{subsec:obstruction},
the obstruction to extending the lift across this cell is the image of $W\wedge V$ in $H_6(M_q,F;\Z)$.
If $W\wedge V\in I$, exactness of \eqref{common:primary-obstruction} makes this obstruction zero and the lift extends.
Conversely, suppose that $h:S^3\times S^3\to F$ extends $s\vee t$. Then $q\circ h$ and $f$ agree on the wedge.
Since $\pi_6(K_{\R,3})=0$, they are homotopic on the six-cell relative to its boundary.
Hence $W\wedge V=q_*h_*[S^3\times S^3]\in I$.
\end{proof}

\begin{lemma}\label{h3:geo:bounded-kernel}
Put $J_q=H_6(Z_q;\Z)$. There is a short exact sequence
\[
 0\longrightarrow J_q\longrightarrow H_6(F;\Z)\xrightarrow{q_*}I\longrightarrow0.
\]
Via $\pi_6F\cong\pi_6Z_q$, the Hurewicz map of $Z_q$ gives a surjection $h:\pi_6F\to J_q$ with $2\ker h=0$.
\end{lemma}
\begin{proof}
Since the base $K_{\R,3}$ is 2-connected and the fiber $Z_q$ is 4-connected, only the two axes of the integral
Serre spectral sequence can be nonzero in total degrees at most seven. Since $H_7(K_{\R,3};\Z)=0$, we have $E^\infty_{0,6}=J_q$. The two
edge homomorphisms in total degree six yield the short exact sequence.

Since $Z_q$ is 4-connected, the low-degree part of Whitehead's exact sequence \cite[\S10 and the end of \S14, p.~81]{JHCWhitehead} reads
\[
 \pi_5Z_q\otimes_{\Z}\Z/2\Z\longrightarrow
 \pi_6Z_q\xrightarrow{\ h\ }H_6(Z_q;\Z)\longrightarrow0.
\]
Identifying $\pi_6Z_q\cong\pi_6F$ gives surjectivity of $h$ and $2\ker h=0$.
\end{proof}

Fix a based sphere map $\sigma_V:S^3\to F$ of volume $V$, and denote the homomorphism induced by point evaluation by
$\ev_*:\pi_3\Map(S^3,F)_{\sigma_V}\to\pi_3F\cong\R$.
\begin{proposition}\label{h3:geo:coherent-lift}
Let $V>0$ and $W\in\R$. Suppose a homomorphism $\gamma:\Q\to H_6(F;\Z)$ satisfies
\[
 q_*\gamma(r)=r(W\wedge V)\qquad(r\in\Q).
\]
Then there is a homomorphism $\beta:\Q\to\pi_3\Map(S^3,F)_{\sigma_V}$ such that $\ev_*\beta(r)=rW$ for all $r\in\Q$.
\end{proposition}
\begin{proof}
Put $A=\pi_3\Map(S^3,F)_{\sigma_V}$ and $K_V=\{W:W\wedge V\in I\}$. Since $\pi_1F=0$, the fiber over the
basepoint of the evaluation map on $\Map(S^3,F)_{\sigma_V}$ is $\Map_*(S^3,F)_{\sigma_V}$.
Translation by loop multiplication identifies the third homotopy group of this fiber with $\pi_6F$.
The homotopy exact sequence, together with $\pi_4F=0$ and Lemma~\ref{h3:geo:primary-obstruction}, now gives
\[
 0\longrightarrow\pi_6F\xrightarrow{j}A\xrightarrow{\ev_*}K_V\longrightarrow0.
\]
Let $\tau:A\to H_6(F;\Z)$ send a class to the image of the fundamental class under its adjoint $S^3\times S^3\to F$. This is additive,
being the composite of the Hurewicz map, cross product with $[S^3]$ in the second factor, and the map
induced by evaluation. Taking the parameter sphere first gives
\[
 q_*\tau(a)=\ev_*(a)\wedge V,\qquad \tau j=h.
\]
The first identity follows by applying \eqref{h3:src:finite-wedge} to the two pure components on $S^3\times S^3$. In the second identity we
regard $J_q\subset H_6(F;\Z)$ as a subgroup. For the translation between components use the pinch map $S^3\to S^3\vee S^3$ whose
projections to both summands have degree $+1$. Send $g$ to the composite of this pinch map with $\sigma_V\vee g$, and
use a based homotopy to identify the image of the constant map with $\sigma_V$. The $\sigma_V$ part of the adjoint
factors through the second sphere and contributes nothing in degree six. The $g$ part factors through
the quotient $S^3\wedge S^3=S^6$ with its product orientation and gives $h$.

Thus $p=(\tau,\ev_*)$ maps onto $\{(c,w)\in H_6(F;\Z)\times K_V:q_*c=w\wedge V\}$. Indeed, choose
$a\in A$ with $\ev_*(a)=w$. Then $c-\tau(a)$ lies in $J_q=\im(\tau j)$, so we may adjust $a$ by an
element of $\im j$ to obtain $p(a)=(c,w)$. Moreover, $\ker p=j(\ker h)$, and hence $2\ker p=0$ by
Lemma~\ref{h3:geo:bounded-kernel}. Define $\beta(r)$ to be twice any preimage of $(\gamma(r/2),rW/2)$.
Any two preimages differ by an element killed by $2$, so $\beta(r)$ is well defined. Adding preimages
proves additivity, and $p\beta(r)=(\gamma(r),rW)$.
\end{proof}

We now take $V=\totalVol$ and write $\sigma_{\totalVol}$ for the reference map $S^3\to F$.
Choose a based map $c_M:M\to S^3$ of degree one by collapsing the complement of an oriented ball disjoint from $a_0$.
Since $F$ is $2$-connected and $\dim M=3$, pulling back $u_F$ gives a bijection $[M,F]\cong H^3(M;\R)$.
Thus $\sigma_{\totalVol}c_M\simeq\sigma$. Precomposition with $c_M$ and change of basepoint along this homotopy
give a homomorphism from $\pi_3\Map(S^3,F)_{\sigma_{\totalVol}}$ to $\pi_3\mathcal S_M$.
Lifting to the universal cover and applying the integral Hurewicz map gives a class in $H_3(\widetilde{\mathcal S}_M;\Z)$.
Composition with $\mathsf T_{\cH}$ now gives a homomorphism
$\Psi:\pi_3\Map(S^3,F)_{\sigma_{\totalVol}}\to H_3(B\cH^\delta;\Z)$.
Evaluation at the basepoint commutes with precomposition by $c_M$, and change of basepoint preserves the volume pairing.
The $\cH$ formula in Proposition~\ref{common:comparison} therefore gives $P_{\cH}\Psi(a)=\ev_*(a)$.

\begin{theorem}\label{h3:geo:integer-additive-section}
There is a noncanonical additive injection $\alpha_{\cH}^{\mathcal P_{\mathrm{div}}}:\mathcal P_{\mathrm{div}}\to H_3(B\cH^\delta;\Z)$. Put $\alpha_{\cG}^{\mathcal P_{\mathrm{div}}}\coloneqq i_*\alpha_{\cH}^{\mathcal P_{\mathrm{div}}}$. Then for both $\cL=\cH,\cG$ we have
\begin{equation}\label{h3:geo:alpha-formula}
 P_{\cL}(\alpha_{\cL}^{\mathcal P_{\mathrm{div}}}(u))=\totalVol\,u\qquad(u\in\mathcal P_{\mathrm{div}})
\end{equation}
and $\alpha_{\cG}^{\mathcal P_{\mathrm{div}}}$ is injective as well.
\end{theorem}
\begin{proof}
Choose a rational basis $\{u_i\}$ of $\mathcal P_{\mathrm{div}}$ and homomorphisms $\gamma_i:\Q\to H_6(F;\Z)$ satisfying $q_*\gamma_i(1)=1\wedge u_i$. Since $q_*\gamma_i(r)=r(1\wedge u_i)$, we have
\[
 q_*\bigl(-(h_{\totalVol})_*\gamma_i(r)\bigr)=r(\totalVol\,u_i\wedge \totalVol).
\]
Proposition~\ref{h3:geo:coherent-lift} gives an additive family $\beta_i$ whose point evaluation has volume $r\,\totalVol\,u_i$. Define $\alpha_{\cH}^{\mathcal P_{\mathrm{div}}}(\sum_i r_i u_i)=\sum_i\Psi\beta_i(r_i)$. Then $P_{\cH}\alpha_{\cH}^{\mathcal P_{\mathrm{div}}}(u)=\totalVol\,u$.
Since $P_{\cG}i_*=P_{\cH}$, the periods on $\cG$ are the same, and $\totalVol>0$ makes both families injective.
\end{proof}

\subsection{Real closedness and ratios of logarithms}
\begin{lemma}\label{h3:src:cubic}
For every $a\in\R$ we have
\[
 \mathfrak r(a)=1\wedge a^3-3a\wedge a^2\in I_{\mathrm{div}}.
\]
\end{lemma}
\begin{proof}
Let $(u,v)\in\Q^3\times\Q^3$ act on $\R^3$ by translation through $u+av$.
This action preserves the foliation by points and the standard transverse volume. Its derivative is the identity,
so the action defines a framed classifying map $B\Q^6\to F$. Denote the degree-one coordinate classes by
$u_i,v_i$. The pullback of the volume class is then $(u_1+av_1)(u_2+av_2)(u_3+av_3)$.
Use the coordinate order $(u_1,u_2,u_3,v_1,v_2,v_3)$ to identify $H_6(B\Q^6;\Z)=\bigwedge^6_\Q\Q^6$ with $\Q$.
In \eqref{h3:src:finite-wedge}, the two pure components contribute $1\wedge a^3$ to the image of $1$,
and each of the three complementary pairs of mixed components contributes $-a\wedge a^2$.
The resulting element belongs to $I_{\mathrm{div}}$ because the classifying map induces a homomorphism
from $\Q=H_6(B\Q^6;\Z)$ to $H_6(F;\Z)$.
\end{proof}

\begin{lemma}\label{h3:src:algebra}
Let $L\subset\R\wedge_{\Q}\R$ be a rational vector subspace invariant under $T_\lambda$ for all $\lambda>0$ and containing every $\mathfrak r(a)$. Write $A_L$
for the quotient $(\R\wedge_{\Q}\R)/L$ and use brackets for the quotient map. Then $A_L$ has a natural real vector space
structure, and $D_L(t)=[1\wedge t]$ is a $\Q$-derivation with values in $A_L$. Moreover,
\[
 [s\wedge t]=sD_L(t)-tD_L(s).
\]
In particular $\{t:1\wedge t\in L\}=\ker D_L$ is a real closed subfield.
\end{lemma}
\begin{proof}
Put $T=\R\otimes_\Q\R$, and write $T^\pm$ for the eigenspaces of the interchange of factors. Identify the exterior square
with $T^-$ by $s\wedge t\mapsto s\otimes t-t\otimes s$. Since $T^+$ is the rational span of the elements $a\otimes a$, invariance under dilations makes $L$
a $T^+$-submodule. Set
\[
 m_a=\tfrac12(a\otimes1+1\otimes a),\qquad d_a=1\otimes a-a\otimes1.
\]
The cubic relation $\mathfrak r(a)\in L$ reads $d_a^3\in L$. Polarization gives $d_ad_bd_c\in L$. We shall use the identities
\[
 s\otimes t-t\otimes s=m_s d_t-m_t d_s,\qquad
 d_{ab}=m_a d_b+m_b d_a,\qquad m_a m_b-m_{ab}=-\tfrac14 d_a d_b.
\]
By the first identity, $T^-$ is generated over $T^+$ by the elements $d_c$, so the right-hand side of the
third identity acts as zero on $T^-/L$. Multiplication by $m_a$ therefore defines the action of the scalar
$a$, and the second identity makes $D_L(a)=[d_a]$ a derivation. The first identity gives $[s\wedge t]=sD_L(t)-tD_L(s)$.
The kernel of a derivation is a subfield. Let $x\in\R$ be algebraic over $\ker D_L$, with minimal polynomial $p$.
Applying $D_L$ to $p(x)=0$ gives $p'(x)D_L(x)=0$. By separability, $p'(x)\ne0$, hence $D_L(x)=0$.
Thus $\ker D_L$ is relatively algebraically closed in $\R$ and is therefore real closed.
\end{proof}

To realize ratios of logarithms, we use the diagonal right action of the real Sol group on its square.
We write the group and its multiplication as
\[
 S=\R^2\rtimes\R,\qquad
 (x,y,z)(x',y',z')=(x+e^zx',y+e^{-z}y',z+z')
\]
and take the left-invariant forms and bi-invariant volume form to be
\[
 \theta_x=e^{-z}dx,\quad\theta_y=e^zdy,\quad\theta_z=dz,
 \qquad\eta=\theta_x\theta_y\theta_z=dx\wedge dy\wedge dz.
\]
Juxtaposition of differential forms denotes exterior product.

\begin{lemma}\label{h3:src:diagonal}
Let $\Lambda$ be a discrete group with two homomorphisms $\rho_1,\rho_2:\Lambda\to S$.
On $S\times S$, the foliation by diagonal right $S$-orbits, equipped with the transverse volume induced by $\eta$,
defines a framed classifying map $f_\rho:B\Lambda\to F$ for the left $\Lambda$-action given by $\rho_1,\rho_2$.
Let $\zeta_i$ and $B_i$ be the classes on $B\Lambda$ determined by the invariant forms $dz_i$ and $dx_i\wedge dy_i$,
respectively. Then
\begin{equation}\label{h3:src:diagonal-volume}
 f_\rho^*u_F=(\zeta_2-\zeta_1)(B_1+B_2).
\end{equation}
\end{lemma}
\begin{proof}
The right diagonal orbits are the fibers of $\tau(g_1,g_2)=g_2g_1^{-1}$. An element $g\in\Lambda$ acts on $p=(g_1,g_2)$ by $gp=(\rho_1(g)g_1,\rho_2(g)g_2)$, and $\tau(gp)=\rho_2(g)\tau(p)\rho_1(g)^{-1}$, so $\tau^*\eta$ is
invariant. Put $\vartheta=(\theta_x,\theta_y,\theta_z)$. Then $\beta=\vartheta_2-\vartheta_1$ is an invariant conormal coframe, $\tau^*\vartheta=\operatorname{Ad}_{g_1}\beta$, and $\det\operatorname{Ad}_{g_1}=1$. Send an arrow $(g,p)$ of the
action groupoid to the germ at $\tau(p)$ of $s\mapsto\rho_2(g)s\rho_1(g)^{-1}$. This defines a continuous functor. The isomorphism $T(S\times S)/\ker d\tau\cong\tau^*TS$
induced by $d\tau$ sends the frame dual to $\beta$ to a volume frame equivariant for the derivative action.
We thus obtain a map from $E\Lambda\times_\Lambda(S\times S)\simeq B\Lambda$ to $F$ \cite[pp.~534--536]{McDuffIII}.

Since the left action is affine, we may send each vertex $g_j$ of a homogeneous bar simplex to $(\rho_1(g_j)0,\rho_2(g_j)0)$ and
extend in barycentric coordinates. Integration of invariant forms over these simplices takes wedge
products to cup products in cohomology \cite[Theorem~2.33]{Dupont}. The structure equations $d\theta_x=\theta_x\theta_z$ and $d\theta_y=-\theta_y\theta_z$ give
\[
 \tau^*\eta=(\theta_{z,2}-\theta_{z,1})
 (\theta_{x,1}\theta_{y,1}+\theta_{x,2}\theta_{y,2})
 +d(\theta_{x,2}\theta_{y,1}-\theta_{x,1}\theta_{y,2}) .
\]
The form inside the last differential is invariant, so its differential represents a coboundary in the bar complex.
This proves \eqref{h3:src:diagonal-volume}.
\end{proof}

For a real algebraic number $a>1$ put $K_a=\Q(a)$, $A_a=K_a^2$, and $\Gamma_a=A_a\rtimes\Z$, with the generator of the quotient acting by $D_a=\operatorname{diag}(a,a^{-1})$.
Write $\ell_a=\log a$. Then $(x,y,n)\mapsto(x,y,n\ell_a)$ defines a homomorphism $\Gamma_a\to S$.
Let $C_a$ be the pullback to $B\Gamma_a$ of the positive generator of $H^1(B\Z;\Z)$ under the quotient map.
Using $\Gamma_a\to S$ as one of the homomorphisms in Lemma~\ref{h3:src:diagonal} gives the degree-one class
$\ell_a C_a$. Let $B_a\in H^2(B\Gamma_a;\R)$ be the class determined by $dx\wedge dy$ in the same factor.
For translations $r=(x,y)$ and $s=(x',y')$, evaluating $B_a$ on the two affine triangles of the bar cycle
$[r|s]-[s|r]$ gives $xy'-yx'$.

\begin{lemma}\label{h3:src:casimir}
There exists a class $c_a\in H_3(B\Gamma_a;\Z)$ such that $\langle B_aC_a,c_a\rangle=1$.
Moreover, the abelian group $H_3(B\Gamma_a;\Z)$ is uniquely divisible and hence is a $\Q$-vector space.
\end{lemma}
\begin{proof}
Taking the directed limit over finitely generated free abelian subgroups gives $H_j(A_a;\Z)=\bigwedge^j_\Q A_a$ for $j>0$. The minimal
polynomials of $a,a^{-1}$ are separable, so $D_a$ is diagonalizable over $\overline{\Q}$. Wedges of eigenvectors also
diagonalize $D_{a*}=\bigwedge^2 D_a$. Hence $H_2(A_a;\Z)$ splits over $\Q$ as the direct sum of the kernel and image of $D_{a*}-1$. Since $D_a$
preserves area, evaluation of $B_a$ vanishes on the image. Project $(1,0)\wedge(0,1)$ onto the kernel along this
splitting and denote the result by $e_a$. Then $\langle B_a,e_a\rangle=1$.

In the mapping torus model $BA_a\to B\Gamma_a\to S^1$ for the quotient $\Gamma_a\to\Z$, orient the base circle so that $C_a$ evaluates
positively, and use the integral Wang sequence to lift $e_a$ to $c_a$. The Wang map followed by inclusion
of the fiber is cap product with $C_a$, so $\langle B_aC_a,c_a\rangle=1$. The same sequence expresses $H_3(B\Gamma_a;\Z)$ as an extension of $\ker(\bigwedge^2_{\Q}D_a-1)$ by
$\operatorname{coker}(\bigwedge^3_{\Q}D_a-1)$. Both groups are uniquely divisible, hence so is the
middle group.
\end{proof}

\begin{proposition}\label{h3:src:algebraic-image}
Take two real algebraic numbers $a,b>1$. Compose the two projections from $\Gamma_a\times\Gamma_b$ with the homomorphisms to $S$
defined above, and let $f:B(\Gamma_a\times\Gamma_b)\to F$ be the map given by Lemma~\ref{h3:src:diagonal} for this pair of homomorphisms. Then
\[
 q_*f_*(c_a\times c_b)=-2\ell_a\wedge\ell_b.
\]
In particular, $1\wedge(\log b/\log a)\in I_{\mathrm{div}}$.
\end{proposition}
\begin{proof}
By \eqref{h3:src:diagonal-volume} we may write $f^*u_F$ as
\[
 -\ell_a B_aC_a\times1+1\times\ell_b B_bC_b
 +C_a\times(-\ell_a B_b)+(\ell_b B_a)\times C_b.
\]
In degree three the signs from interchanging the cup factors and the exterior coefficient factors cancel.
Thus $\smile_\wedge$ is symmetric in cohomology. In \eqref{h3:src:finite-wedge}, two pairs of terms contribute
to bidegree $(3,3)$, and the factor $\tfrac12$ allows us to count each pair once.
The two pure terms contribute $(-\ell_a)\wedge\ell_b$ by Lemma~\ref{h3:src:casimir}.
For the mixed terms, $C_a$ and $C_b$ are integral and the sign from interchanging factors is
$(-1)^{2\cdot2}=1$. Their contribution is therefore
\[
 \langle-\ell_a B_b C_b,c_b\rangle\wedge\langle\ell_b C_a B_a,c_a\rangle
 =(-\ell_a)\wedge\ell_b.
\]

The rational vector space structure on $H_3(B\Gamma_a;\Z)$ makes $\gamma(r)=(h_{1/\ell_a})_*f_*\bigl((-\tfrac r2c_a)\times c_b\bigr)$ a homomorphism $\Q\to H_6(F;\Z)$. Since $q_*\gamma(r)=r\bigl(1\wedge(\ell_b/\ell_a)\bigr)$, the claimed inclusion
follows.
\end{proof}

\begin{theorem}\label{h3:src:logfield}
The sets $\mathcal P$ and $\mathcal P_{\mathrm{div}}$ are real closed subfields of $\R$. For the field $E_{\log}$ defined in Section~\ref{sec:intro}, we have
$E_{\log}\subseteq\mathcal P_{\mathrm{div}}$.
\end{theorem}
\begin{proof}
By Lemma~\ref{h3:src:cubic} and $I_{\mathrm{div}}\subseteq I$, Lemma~\ref{h3:src:algebra} applies to both $I$ and $I_{\mathrm{div}}$. Thus $\mathcal P$ and $\mathcal P_{\mathrm{div}}$ are real closed subfields.
By Proposition~\ref{h3:src:algebraic-image}, $\mathcal P_{\mathrm{div}}$ contains $\log b/\log a$ for all real algebraic $a,b>1$.
The defining minimality of $E_{\log}$ gives $E_{\log}\subseteq\mathcal P_{\mathrm{div}}$.
\end{proof}

\begin{proof}[Proof of Theorem~\ref{thm:threeD-E}]
By Theorem~\ref{h3:src:logfield} it suffices to restrict the additive family of Theorem~\ref{h3:geo:integer-additive-section} to $E_{\log}$. The norm bound follows
from \eqref{h3:geo:alpha-formula} and Corollary~\ref{common:lone}. If a family from Theorem~\ref{thm:odd} has already been fixed, choose a rational basis
of $\mathcal P_{\mathrm{div}}$ containing $1$ and use the prescribed $\alpha_{\cH}^{\Q}$ on $\Q\cdot1$. This extends the given family.
\end{proof}

For a closed real hyperbolic three-manifold with its hyperbolic volume form, \eqref{common:bm} now shows that
$\operatorname{comp}(\mathcal B_M)$ evaluates to $-\totalVol^2u$ on $J_*\alpha_{\cG}^{E_{\log}}(u)$.

On both $\cH$ and $\cG$ the additive family of Theorem~\ref{thm:threeD-E} realizes the algebraic number $\sqrt2$ and the
transcendental number $\log 3/\log 2$ \cite[p.~24, discussion after Theorem~2]{Waldschmidt1978} as normalized integral periods.

\subsection{Normalized periods of the flux kernel}\label{h3:periods-section}

\begin{lemma}\label{h3:geo:diffeo-contractible}
Under the hypotheses of Theorem~\ref{thm:threeD-E} both $\cG$ and $\cH$ are contractible. Moreover, the inclusion $\Diff_v(M)\hookrightarrow\Diff^+(M)$
induces an isomorphism of component groups, where $\Diff^+(M)$ denotes the group of orientation-preserving
diffeomorphisms.
\end{lemma}
\begin{proof}
For closed oriented aspherical three-manifolds other than non-Haken infranilmanifolds, McCullough and
Soma \cite[pp.~328--329, Theorem~1, Corollary~1]{McCulloughSoma} give $\Diff_0(M)\simeq(S^1)^k$, where $k$ is the rank of the center of the fundamental group. The excluded manifolds
have $k=1$, so the result applies in the centerless case and $\Diff_0(M)$ is contractible. Let $\mathcal V_{\totalVol}$ denote the space
of positive smooth volume forms of total volume $\totalVol$. The parametric Moser argument identifies $\Diff_v(M)\subset\Diff^+(M)$ with
the slice $\Diff_v(M)\times\{v\}$ of the product $\Diff_v(M)\times\mathcal V_{\totalVol}$ \cite[Part~I, Theorem~2.7.2, p.~167]{MEF}. Since $\mathcal V_{\totalVol}$ is convex, $\Diff_v(M)$ is a deformation retract of $\Diff^+(M)$.
Restriction to the identity components gives contractibility of $\cG$, and the flux decomposition of
Section~\ref{subsec:flat} gives contractibility of $\cH$.
\end{proof}

\begin{proposition}\label{h3:period-characterization}
For $(M,v)$ satisfying the hypotheses of Theorem~\ref{thm:threeD-E}, put $\widehat P_{\cH}\coloneqq P_{\cH}/\totalVol$. Then
\[
 \im\widehat P_{\cH}=\mathcal P,\qquad
 \bigl\{\widehat P_{\cH}(a(1)):a\in\Hom_{\Z}(\Q,H_3(B\cH^\delta;\Z))\bigr\}
 =\mathcal P_{\mathrm{div}}.
\]
Both sets are real closed subfields independent of $M$ and $v$.
\end{proposition}
\begin{proof}
By Lemma~\ref{h3:geo:diffeo-contractible}, $B_{\mathrm{top}}\cH$ is weakly contractible, so $p_{\cH}$ is a weak equivalence. The integral comparison $\mathsf T_{\cH}$ of
Section~\ref{subsec:flat} is therefore an isomorphism. Using the evaluation map $e:\mathcal S_M\times M\to F$ in the chosen volume framing and
the universal covering projection $\rho$, define a homomorphism on integral homology by
\[
 \mathcal T(z)=\bigl(e\circ(\rho\times\mathrm{id}_M)\bigr)_*
               \bigl((\mathsf T_{\cH}^{-1}z)\times[M]\bigr)\in H_6(F;\Z).
\]

Since $\widetilde{\mathcal S}_M$ is simply connected, its degree-three Hurewicz map is surjective
\cite[supplementary Chapter~5, \S5.A, p.~605]{Hatcher}. Choose a map $S^3\to\widetilde{\mathcal S}_M$ representing $\mathsf T_{\cH}^{-1}z$.
The associated evaluation map $S^3\times M\to F$ pulls $u_F$ back to a class with no mixed components, since
$H^1(S^3;\R)=H^2(S^3;\R)=0$. Its restriction to the parameter sphere has period $P_{\cH}(z)$ by Proposition~\ref{common:comparison}.
Its restriction to $M$ has period $\totalVol$ since $\sigma^*u_F=[v]$. Applying \eqref{h3:src:finite-wedge} and then
rescaling the transverse volume gives
\[
 q_*\mathcal T(z)=P_{\cH}(z)\wedge \totalVol,\qquad
 q_*\bigl(-(h_{1/\totalVol})_*\mathcal T(z)\bigr)=1\wedge\widehat P_{\cH}(z).
\]
Thus $\widehat P_{\cH}(z)\in\mathcal P$. For an additive family $a$, the composite $-(h_{1/\totalVol})_*\mathcal T\circ a$ gives coherent integral lifts of the rational
multiples of $1\wedge\widehat P_{\cH}(a(1))$.

Conversely, put $V=\totalVol$. If $u\in\mathcal P$, then $Vu\wedge V=-T_V(1\wedge u)\in I$. By Lemma~\ref{h3:geo:primary-obstruction} the wedge of a sphere map of volume $Vu$ and $\sigma_V$ extends
over the product. Let $b$ be the homotopy class adjoint to this extension. Then $P_{\cH}(\Psi(b))=\ev_*(b)=Vu$ gives $\widehat P_{\cH}(\Psi(b))=u$. If $u\in\mathcal P_{\mathrm{div}}$,
the map $t\mapsto\alpha_{\cH}^{\mathcal P_{\mathrm{div}}}(tu)$ of Theorem~\ref{h3:geo:integer-additive-section} gives the required additive family. Real closedness follows from
Theorem~\ref{h3:src:logfield}. Independence of $(M,v)$ follows from the definitions of $I$ and $I_{\mathrm{div}}$.
\end{proof}

If $b_1(M)=0$, then $\cG=\cH$ by Section~\ref{sec:intro}, so both identities in
Proposition~\ref{h3:period-characterization} hold for the identity component as well.
If $b_1(M)>0$, we can apply Proposition~\ref{oddreal:main} with $p=1$, since all positive-degree rational
Pontryagin classes vanish in dimension three. Thus $\im(P_{\cG}/\totalVol)=\R$.

\section{Real hyperbolic manifolds and transfer to the full group}\label{sec:applications}

We first compare the period groups of the identity component and the flux kernel for closed real
hyperbolic manifolds. We then use transfer to pass to the full volume-preserving diffeomorphism group,
assuming its outer action on the fundamental group has finite image. For a closed real hyperbolic
three-manifold, we evaluate the transferred class on the coherent family of
Section~\ref{sec:threeD-refinement}.

\subsection{Period groups of closed real hyperbolic manifolds}

Throughout this subsection $M$ is a closed connected oriented real hyperbolic manifold of dimension $n\ge3$,
and $v$ is any positive smooth volume form. By Section~\ref{subsec:bm-proof}, $M$ satisfies the geometric
hypotheses of Theorems~\ref{thm:odd} and~\ref{thm:even}. In odd dimensions, substituting
$\|M\|=V_{\rm hyp}/v_n$ into \eqref{odd:eq:norm} gives the lower bound $V_{\rm hyp}|t|/v_n$ for both groups.
This bound is independent of $v$.

\begin{corollary}[Even-dimensional hyperbolic manifolds]\label{apps:four}
Let $M^n$ be a closed connected oriented real hyperbolic manifold of dimension $n=2m\ge4$, and let $v$ be any
positive smooth volume form. Then
\[
 \im P_{\cH}=
 \begin{cases}
 \R & (b_{[2,m]}(M)>0),\\
 \{0\} & (b_{[2,m]}(M)=0),
 \end{cases}
 \qquad
 \im P_{\cG}=
 \begin{cases}
 \im P_{\cH} & (n\equiv0\bmod4),\\
 \R & (n\equiv2\bmod4).
 \end{cases}.
\]
Whenever the integral period group is zero, the corresponding top volume class vanishes as well.
\end{corollary}
\begin{proof}
Use the Gauss--Bonnet sign from Section~\ref{subsec:bm-proof}. If $b_{[2,m]}(M)=0$, Poincar\'e duality gives $\chi(M)=2-2b_1(M)$. Hence $b_1=0$ when $m$ is
even, and $b_1\ge2$ when it is odd. The result follows from $b_{[1,m]}(M)=b_{[2,m]}(M)+b_1(M)$ and Theorem~\ref{thm:even}.
\end{proof}

In particular, in dimension four every real number occurs as a period on both groups if and only if
$b_2(M)>0$.

\subsection{Transfer to the full group}\label{apps:full-group}

Return to the general setting of Section~\ref{sec:intro}, with $n\ge3$ and $D=\Diff_v(M)$. Write $K_{\mathrm{out}}$ for the kernel of the outer
action $D\to\operatorname{Out}(\Gamma)$. Since $M$ is a $K(\Gamma,1)$, this kernel consists precisely of the elements homotopic to the identity
as maps $M\to M$ \cite[Proposition~1B.9]{Hatcher}. Assume the outer action has finite image and put $d=[D:K_{\mathrm{out}}]$. Write $j:\cG\hookrightarrow D$ for the inclusion. The
construction in Section~\ref{subsec:flat} gives a top volume class $\xi_{K_{\mathrm{out}}}$ on $K_{\mathrm{out}}$.
Let $\operatorname{cor}_{K_{\mathrm{out}}}^D$ denote the cohomological transfer for the corresponding regular
$d$-sheeted covering, with the convention that one sums over the sheets. We normalize this transfer to define
$\xi_D\coloneqq d^{-1}\operatorname{cor}_{K_{\mathrm{out}}}^D\xi_{K_{\mathrm{out}}}$ \cite[\S3.G, pp.~321--322]{Hatcher}.

\begin{proposition}[Detecting integral periods on the full group]\label{prop:full-transfer}
We have $\xi_D|_{K_{\mathrm{out}}}=\xi_{K_{\mathrm{out}}}$. Define $P_D(y)\coloneqq\langle\xi_D,y\rangle$. For all $z\in H_n(B\cG^\delta;\Z)$ and $y\in H_n(BD^\delta;\Z)$ we have
\[
 P_D(j_*z)=P_{\cG}(z),\qquad
 \|y_\R\|_1\ge\frac{\|M\|}{\totalVol}|P_D(y)|.
\]
Under the hypotheses of Theorem~\ref{thm:odd}, $j_*\alpha_{\cG}^{\Q}$ is injective and its image is a direct summand of $H_n(BD^\delta;\Z)$.
Under the hypotheses of Theorem~\ref{thm:threeD-E}, $j_*\alpha_{\cG}^{E_{\log}}$ likewise embeds $E_{\log}$ as a direct summand of $H_3(BD^\delta;\Z)$.
In even dimension $n=2m\ge4$, the class $\xi_D$ vanishes if $b_{[1,m]}(M)=0$.
For a rationally Pontryagin-trivial $M$ of even dimension $n=2m\ge4$, we have $\im P_D=\R$ if and only if $b_{[1,m]}(M)>0$.
Under the hypotheses of Proposition~\ref{scope:bott}, one also has $\xi_D=0$.
\end{proposition}
\begin{proof}
The flat bundle map $\widehat c_h$ induced by conjugation by $h\in D$ preserves the transverse volume class and the
fiber orientation. Since the restriction of $h^{-1}\circ r_{K_{\mathrm{out}}}\circ\widehat c_h$ to the reference fiber is
homotopic to the identity, uniqueness of the projection gives $r_{K_{\mathrm{out}}}\circ\widehat
c_h\simeq h\circ r_{K_{\mathrm{out}}}$. Naturality of fiber integration and $h^*[v]=[v]$ show that $\xi_{K_{\mathrm{out}}}$ is invariant under conjugation.
Transfer followed by restriction is the sum of the conjugation actions of $D/K_{\mathrm{out}}$. Normalizing by $d^{-1}$
therefore gives $\xi_D|_{K_{\mathrm{out}}}=\xi_{K_{\mathrm{out}}}$. We also have $\xi_{K_{\mathrm{out}}}|_{\cG}=\xi_{M,n}$.

Write $t$ for integral homology transfer for this regular cover. Then $\langle\xi_{K_{\mathrm{out}}},ty\rangle=dP_D(y)$. Summing the $d$ lifts of each
simplex gives $\|(ty)_\R\|_1\le d\|y_\R\|_1$, so the lower bound follows from Corollary~\ref{common:lone} applied to $K_{\mathrm{out}}$. Period evaluation shows
that the integral families remain injective, and their divisible images are direct summands.

In even dimension $n=2m$, if $b_{[1,m]}(M)=0$ then Poincar\'e duality forces all intermediate Betti numbers to vanish.
Applying \eqref{even:quadraticformula} to $\cL=K_{\mathrm{out}}$ gives $\xi_{K_{\mathrm{out}}}=0$, and hence $\xi_D=0$.
If $b_{[1,m]}(M)>0$ and $M$ is rationally Pontryagin-trivial, Theorem~\ref{thm:even} gives surjectivity of $P_{\cG}$.
Since $P_D\circ j_*=P_{\cG}$, $P_D$ is surjective as well. Finally, under the hypotheses of
Proposition~\ref{scope:bott}, that proposition gives $\xi_{K_{\mathrm{out}}}=0$, so again $\xi_D=0$.
\end{proof}

If $M$ has odd dimension and is rationally Pontryagin-trivial, and $H^p(M;\Q)\ne0$ for some $1\le p\le n-1$, then Proposition~\ref{oddreal:main}
and $P_D\circ j_*=P_{\cG}$ imply that $P_D$ is surjective. For closed real hyperbolic manifolds the outer action has finite
image by Mostow rigidity \cite{Mostow}.

\begin{corollary}[The full group of a closed real hyperbolic three-manifold]\label{apps:hyperbolic-three-transfer}
Let $M$ be a closed connected oriented real hyperbolic three-manifold with any positive smooth volume
form $v$. Denote by $\operatorname{Isom}^+(M)$ the group of orientation-preserving isometries of its hyperbolic metric. Then
$K_{\mathrm{out}}=\cG$ and $d=[D:\cG]=|\operatorname{Isom}^+(M)|$. The map $j_*\alpha_{\cG}^{E_{\log}}$ is injective
and its image is a direct summand of $H_3(BD^\delta;\Z)$. For the unnormalized transfer we have
\begin{equation}\label{apps:sum-transfer}
 \bigl\langle\operatorname{cor}_{\cG}^D\xi_{M,3},j_*\alpha_{\cG}^{E_{\log}}(u)\bigr\rangle=d\,\totalVol\,u,
 \qquad
 \|(j_*\alpha_{\cG}^{E_{\log}}(u))_\R\|_1\ge\|M\|\,|u|
 \qquad(u\in E_{\log}).
\end{equation}
In particular, evaluation at $u=1$ shows that the transferred ordinary class is nonzero.
\end{corollary}
\begin{proof}
By Lemma~\ref{h3:geo:diffeo-contractible} the inclusion $D\hookrightarrow\Diff^+(M)$ induces an isomorphism on component groups. By Gabai \cite[Theorem~1.1]{Gabai} the inclusion of
the isometry group induces an isomorphism $\operatorname{Isom}^+(M)\cong\pi_0\Diff^+(M)$. Mostow rigidity makes the map from this isometry group
to $\operatorname{Out}(\pi_1M)$ injective. Hence $K_{\mathrm{out}}=\cG$ and $d=|\operatorname{Isom}^+(M)|$. Proposition~\ref{prop:full-transfer} and Theorem~\ref{thm:threeD-E} together with $\operatorname{cor}_{\cG}^D\xi_{M,3}=d\,\xi_D$ give \eqref{apps:sum-transfer}.
\end{proof}

\section{Examples of the period criteria}\label{sec:scope}

We first apply Theorem~\ref{thm:even} to a three-torus bundle over the circle for which the identity component
and the flux kernel have different integral period groups. We then apply Proposition~\ref{oddreal:main} to compact
quotients of the symmetric spaces $\SL_k(\R)/\SO(k)$ for $k\ge4$ and obtain all real periods on both groups. Finally,
closed complex hyperbolic manifolds of complex dimension at least two illustrate the vanishing
criterion in Proposition~\ref{scope:bott}.

\begin{example}\label{scope:torusmapping}
Let $M_A$ be the mapping torus of the automorphism of $\mathbb T^3$ induced by $A$, where
\[
 A=\begin{pmatrix}0&0&1\\1&0&1\\0&1&0\end{pmatrix}\in\SL_3(\Z),
 \qquad
 M_A=(\mathbb T^3\times[0,1])/((x,1)\sim(Ax,0)).
\]
For every positive volume form on this four-manifold we have $\im P_{\cG}=\R$, $\xi_{M_A,4}|_{\cH}=0$, and $\im P_{\cH}=\{0\}$. Moreover, every real
period of $P_{\cG}$ is realized by an integral class whose real image has zero $\ell^1$-seminorm. Consequently $\xi_{M_A,4}$
has no bounded cocycle representative.
\end{example}
\begin{proof}
The characteristic polynomial of $A$ is $t^3-t-1$. It has discriminant $-23$, one real root greater than $1$,
and a nonreal conjugate pair. We obtain a real logarithm $B$ of $A^{-1}$ by taking logarithms of its
positive real eigenvalue and of the real two-dimensional block corresponding to the conjugate pair.
Thus $\exp B=A^{-1}$, and $e^{\operatorname{tr}B}=\det A^{-1}=1$ gives $\operatorname{tr}B=0$.
Put $S=\R^3\rtimes_{\exp(tB)}\R$ and $\Lambda=\Z^3\rtimes_{A^{-1}}\Z$. Left translation by $(0,1)\in\Lambda$ takes $(Ax,0)$ to $(x,1)$, so $M_A=\Lambda\backslash S$. The contractible group $S$ is the
universal cover. A left-invariant frame in the $\R^3$ directions and a complementary left-invariant
vector field descend to the quotient, trivializing both the vertical tangent bundle and $TM_A$. A
central element $(w,k)$ of $\pi_1 M_A=\Lambda$ satisfies $A^k=I_3$ and $Aw=w$. The real root above and $\det(A-I_3)=1$ imply $k=0$ and $w=0$. The
abelianization is $\Z\oplus\Z^3/(A^{-1}-I_3)\Z^3\cong\Z$, so $b_1=1$. Since this is a bundle over the circle, $\chi(M_A)=0$. Poincar\'e duality then
gives $b_2=0$. By Theorem~\ref{thm:even}, $\xi_{M_A,4}|_{\cH}=0$ and $\im P_{\cH}=\{0\}$.

Since $\operatorname{tr}B=0$, the form $dx_1\wedge dx_2\wedge dx_3\wedge dt$ is bi-invariant. It
descends to $M_A$, where we rescale it to total volume $\totalVol$. By Moser's theorem \cite[Part~I,
Theorem~2.7.1]{MEF}, a diffeomorphism carries this form to the given one. Naturality of the top volume
class under this identification allows us to use the rescaled form.

For $a\in\R$ define $m_a:S\to S$ by $m_a(x,t)=(ax,t)$. This is a homomorphism. Define $i_a:\Lambda\to\cG$ by $i_a(\lambda)([s])=[s\,m_a(\lambda)^{-1}]$. Right translations preserve the
bi-invariant volume and belong to the identity component because $S$ is connected. Using the inverse
makes right translation compatible with the usual composition of maps. In the flat bundle pulled back
by $i_a$, an element $\lambda$ takes $(s,x)\in S\times M_A$ to $(\lambda s,x\,m_a(\lambda)^{-1})$. Thus
\[
 [s,x]\longmapsto ([s],x\,m_a(s))
\]
is a bundle trivialization. Its second component restricts to the identity on the reference fiber and hence
is homotopic to $r$. The section that is constant with value $\Lambda e$ in this trivialization is
$[s]\mapsto[s,\Lambda m_a(s)^{-1}]$ in the original coordinates. Pulling the transverse volume class back
along this section gives $i_a^*\xi_{M_A,4}$. Inversion on $S$ multiplies the bi-invariant volume form by
$(-1)^4=1$, and pullback by $m_a$ multiplies it by $a^3$. Using $M_A$ as a model for $B\Lambda$, we obtain
$P_{\cG}((i_a)_*[M_A])=a^3\totalVol$. These periods range over $\R$ as $a$ varies.

The map $[x,t]\mapsto[2x,t]$ is an eight-sheeted self-cover of $M_A$, so $\|M_A\|=0$. Hence every class $((i_a)_*[M_A])_\R$ has vanishing
$\ell^1$-seminorm. The period at $a=1$ is nonzero. A bounded cocycle evaluates to zero on classes with
vanishing real $\ell^1$-seminorm, so $\xi_{M_A,4}$ has no bounded representative.
\end{proof}

\begin{corollary}\label{scope:slfive}
Suppose $k\ge4$ and let $\Lambda<\SL_k(\R)$ be a torsion-free cocompact lattice. For every positive smooth volume form on
$M=\Lambda\backslash\SL_k(\R)/\SO(k)$ we have $\im P_{\cG}=\im P_{\cH}=\R$.
\end{corollary}
\begin{proof}
The space $\SL_k(\R)/\SO(k)$ of positive-definite symmetric matrices of determinant one is contractible. Its
dimension is $n=k(k+1)/2-1\ge9$, with $n=14$ when $k=5$. Each stabilizer is the intersection of $\Lambda$ with a conjugate of $\SO(k)$
and is finite. Hence the torsion-free group $\Lambda$ acts freely. The connected group $\SL_k(\R)$ preserves
orientation, so $M$ is a closed oriented aspherical manifold with $\pi_1M\cong\Lambda$. By Borel density \cite{BorelDensity}, the
center $Z(\Lambda)$ lies in the finite center of $\SL_k(\R)$. Since $\Lambda$ is torsion-free, $Z(\Lambda)=1$. By Tshishiku \cite[Theorem~6.2]{Tshishiku}, $M$ is
rationally Pontryagin-trivial. Franke's calculation \cite[preprint, \S7.1, Proposition~1, p.~36]{Franke} shows that for either parity of $k$, the
complex cohomology of the compact dual $\operatorname{SU}(k)/\SO(k)$ has a generator in degree $5$. Thus Matsushima's
inequality \cite[pp.~1--2]{Matsushima} gives $b_5(M)>0$. Since $2\le5\le n-2$, Proposition~\ref{oddreal:main} with $p=5$ gives the result for both groups.
\end{proof}

Finally, let $M=\Lambda\backslash\mathbb H_{\mathbb C}^m$ be a closed complex hyperbolic manifold, where
$m\ge2$ and $\Lambda<\operatorname{PU}(m,1)$ is a torsion-free cocompact lattice. Then $M$ is aspherical, and
its K\"ahler class gives $b_2(M)>0$. Borel density and the triviality of the center of $\operatorname{PU}(m,1)$
give $Z(\Lambda)=1$. Let $h$ denote the hyperplane class of the compact dual $\mathbb{CP}^m$.
The Euler sequence on $\mathbb{CP}^m$ gives $p_1(T\mathbb{CP}^m)=(m+1)h^2\ne0$.
By proportionality with the compact dual \cite[p.~2707]{Tshishiku}, $p_1(TM)\ne0$.
Proposition~\ref{scope:bott} therefore implies that the top volume class and its integral period group vanish on both groups.

\appendix
\section{Comparison with the local point-evaluation cocycle}\label{app:local}

Let $M$ be a closed connected oriented real hyperbolic three-manifold with hyperbolic volume form $v$ and
total volume $\totalVol=\int_Mv$. Let $\cG$ be the identity component of its volume-preserving diffeomorphism group.
We modify the geodesic volume cocycle $F_{M,3}$ of Proposition~\ref{common:hyperbolic} by a bounded coboundary
so that it agrees near the identity with a prescribed local point-evaluation cocycle.
Write $v_3$ for the supremum of the volumes of geodesic tetrahedra in $\mathbb H^3$.

For a symmetric open neighborhood $U$ of the identity $e\in \cG$ put $S_{k,U}=\{(g_0,\ldots,g_k):g_i^{-1}g_j\in U\ (0\le i,j\le k)\}$. Choose $U$ sufficiently small. By
Nosaka's construction \cite[Proposition~3.2]{Nosaka}, for $0\le k\le3$ there are smooth left-equivariant maps $S_{k,U}\times\Delta^k\to \cG$, $(\mathbf g,t)\mapsto T_k(\mathbf g)(t)$, with vertices $g_i$ and
with restrictions to faces given by the corresponding $T_{k-1}$. Fix $a_0\in M$, write $\ev_{a_0}(g)=g(a_0)$, and put
\[
 \phi_T(\mathbf g)\coloneqq\int_{\Delta^3}T_3(\mathbf g)^*\ev_{a_0}^*v.
\]
Two cochains have the same germ if they agree on an open neighborhood of the identity tuple. Local
smoothness refers to the Fr\'echet--Lie group structure on $\cG$.

\begin{proposition}[A bounded extension with prescribed germ]\label{apps:localextension}
For every $\varepsilon>0$ there is a globally smooth bounded homogeneous three-cocycle $\widehat\Phi_\varepsilon:\cG^4\to\R$ with the same germ as $\phi_T$
such that
\[
 [\widehat\Phi_\varepsilon]_b=[F_{M,3}]_b,\qquad
 [\widehat\Phi_\varepsilon]=\xi_{M,3},\qquad
 \|\widehat\Phi_\varepsilon\|_\infty\le v_3,\qquad
 \|\widehat\Phi_\varepsilon-F_{M,3}\|_\infty\le\varepsilon.
\]
\end{proposition}
\begin{proof}
Put $\bar v=v/\totalVol$ and write $\ev_y(g)=g(y)$. Averaging the point evaluation of $T_3$ over $y$ gives the three-cochain
\[
 \bar a(\mathbf g)=\int_M\left(\int_{\Delta^3}(\ev_y\circ T_3(\mathbf g))^*v\right)\bar v(y).
\]
For both prisms below, give $I\times\Delta^k$ the product orientation with $I=[0,1]$ first.
Since $\partial(I\times\Delta^3)=(\{1\}-\{0\})\times\Delta^3-I\times\partial\Delta^3$, Stokes' theorem
shows that the coboundary of the two-cochain obtained by integration over the prisms $I\times\Delta^2$
is the integral at $s=1$ minus the integral at $s=0$, provided the prisms agree on faces.

For the first prism, vary the point of evaluation. Let $w(z)$ be the coordinate indexed by the identity of
$\Gamma=\pi_1M$ in the probability vector $p(z)$ constructed in the proof of Proposition~\ref{common:continuous}.
Then $w$ is smooth, nonnegative, and compactly supported, with $\sum_{\gamma\in\Gamma}w(\gamma z)=1$.
Consequently $q_M:\mathbb H^3\to M$ sends the measure $d\mu(z)=w(z)\,q_M^*\bar v(z)$ to $\bar v$.
Fix a lift of $a_0$ and project the geodesic segment from this lift to $z$ to obtain a path $\gamma_z:I\to M$.
The maps $C^{(k)}_{z,\mathbf g}(s,t)=T_k(\mathbf g)(t)(\gamma_z(s))$ agree on faces.
Their restrictions to $s=0$ and $s=1$ are the simplices obtained by evaluation at $a_0$ and $q_M(z)$, respectively. Hence
\begin{equation}\label{apps:pointmotion}
 b_{\rm mov}(\mathbf g)=\int_{\mathbb H^3}\!\left(\int_{I\times\Delta^2}(C^{(2)}_{z,\mathbf g})^*v\right)d\mu(z),
 \qquad \delta b_{\rm mov}=\bar a-\phi_T .
\end{equation}

For the second prism, straighten the simplices obtained by evaluating the inverse diffeomorphisms.
Set $I_y(g)=g^{-1}(y)$. For $X\in T_g\cG$, viewed as a vector field along $g$, and $x=g^{-1}(y)$, we have $D_gI_y(X)=-(D_xg)^{-1}X(x)$.
To compute $I_y^*v$, apply this formula to each of the three arguments of $v$. Using $g^*v=v$ and changing variables
by $y=g(x)$ gives $\int_M I_y^*v\,\bar v(y)=-\int_M\ev_y^*v\,\bar v(y)$.
Choose a lift $\widetilde y$ of $y$, put $a(t)=\widetilde{T_k(\mathbf g)(t)}^{-1}\widetilde y$, and let $c(t)$
be the straight simplex of Proposition~\ref{common:hyperbolic} with the same vertices as $a(t)$.
At each $t$, project the geodesic segment from $a(t)$ to $c(t)$ to $M$. The resulting prism
$H^{(k)}_{y,\mathbf g}(s,t)$ is compatible with the face maps and, by deck equivariance, independent of the lift $\widetilde y$.
The averaged integrals at its two ends are $-\bar a$ and $A_{M,3}/\totalVol=-F_{M,3}$,
respectively, by the preceding computation and \eqref{common:average}--\eqref{common:normalized}. Thus
\begin{equation}\label{apps:straightprism}
 b_{\rm str}(\mathbf g)=\int_M\!\left(\int_{I\times\Delta^2}(H^{(2)}_{y,\mathbf g})^*v\right)\bar v(y),
 \qquad \delta b_{\rm str}=\bar a-F_{M,3}.
\end{equation}
Left equivariance of $T_k$ and invariance of $v$ under $\cG$ show that $b_{\rm mov}$ is homogeneous.
The identity $I_y(hg)=I_{h^{-1}y}(g)$ and a change of variables give homogeneity of $b_{\rm str}$.
Both cochains are locally smooth, since the integrands depend smoothly on the parameters and the
integrals are over compact domains or have compact support.

By \eqref{apps:pointmotion} and \eqref{apps:straightprism} the difference $b=b_{\rm str}-b_{\rm mov}$ satisfies $\delta b=\phi_T-F_{M,3}$ locally. Since a constant two-cochain has zero
coboundary, subtract $b(e,e,e)$ to arrange $b(e,e,e)=0$.

By Gromov's geodesic cone inequality \cite[\S1.2, p.~19]{GromovVolume}, the absolute volume of a
geodesic tetrahedron in $\mathbb H^3$ is at most half the area of any face. Choose $\rho>0$ so that a
hyperbolic disk of radius $\rho$ has area less than $v_3$. If two vertices of a face are within $\rho$
of the third, the face lies in such a disk by convexity. The absolute volume of the tetrahedron is
then at most $v_3/2$, wherever the opposite vertex lies.

Choose an open neighborhood $O\subset\cG^2$ of $(e,e)$ so that $(h_1,h_2)\in O$ implies $(e,h_1,h_2)\in S_{2,U}$ and
\[
 |b(e,h_1,h_2)|<\tfrac14\min(\varepsilon,v_3/2),\qquad
 \sup_{z\in\mathbb H^3}d(\widetilde h_i^{-1}z,z)<\rho\quad(i=1,2).
\]
For the last condition it suffices to work on a single compact fundamental domain because the
canonical lifts commute with deck transformations. Choose a smooth function $\chi:\cG^2\to[0,1]$ equal to $1$ near $(e,e)$,
with support in $O$. To construct one, fix a smooth embedding $\iota:M\hookrightarrow\R^N$ and compose $\sum_{i=1}^2\|\iota\circ h_i-\iota\|_{H^s}^2$ with a one-variable
cutoff. For sufficiently large $s$ the Sobolev embedding theorem allows the support to be taken in $O$.
Let $b_0$ be the extension by zero of $\chi(h_1,h_2)b(e,h_1,h_2)$ outside $O$ and put $\widetilde b(g_0,g_1,g_2)=b_0(g_0^{-1}g_1,g_0^{-1}g_2)$. This is a globally smooth bounded
homogeneous two-cochain with the same germ as $b$, and it satisfies
\[
 \|\delta\widetilde b\|_\infty\le4\|\widetilde b\|_\infty\le\min(\varepsilon,v_3/2).
\]

If $\delta\widetilde b(\mathbf g)\ne0$, then $\widetilde b(g_i,g_j,g_k)\ne0$ on some face, so
$(g_i^{-1}g_j,g_i^{-1}g_k)\in O$. For every $\widetilde x$, the vertices $\widetilde g_j^{-1}\widetilde x$ and
$\widetilde g_k^{-1}\widetilde x$ in \eqref{common:average} are within distance $\rho$ of $\widetilde g_i^{-1}\widetilde x$.
Taking this face as the base and averaging the volume estimate against $\bar v$ gives
$|F_{M,3}(\mathbf g)|\le v_3/2$. Put $\widehat\Phi_\varepsilon=F_{M,3}+\delta\widetilde b$.
Where the coboundary is nonzero, each summand has absolute value at most $v_3/2$.
Elsewhere \eqref{common:bounded} gives $|\widehat\Phi_\varepsilon|\le v_3$.
The bound on $\delta\widetilde b$ also gives $\|\widehat\Phi_\varepsilon-F_{M,3}\|_\infty\le\varepsilon$.
Near the identity tuple, $\widehat\Phi_\varepsilon=F_{M,3}+\delta b=\phi_T$.
Both $F_{M,3}$ and $\delta\widetilde b$ are globally smooth, so $\widehat\Phi_\varepsilon$ is globally smooth as well.
Since $\widetilde b$ is bounded, the two cocycles
represent the same bounded cohomology class, whose image in ordinary cohomology is $\xi_{M,3}$ by \eqref{common:bounded}.
\end{proof}

\begin{lemma}\label{app:germ}
Let $L$ be a topological group whose underlying space is contractible. Two real-valued homogeneous
group $q$-cocycles with the same germ represent the same class in $H^q(BL^\delta;\R)$.
\end{lemma}
\begin{proof}
Let $c$ be the difference of the cocycles. Choose an open neighborhood $U$ of the identity such that $c$ vanishes
on $U^{q+1}$. By homogeneity, $c$ vanishes on $(gU)^{q+1}$ for every $g\in L$.
Consider the singular chain complex generated by simplices whose images lie in some $gU$.
By the small chains theorem \cite[Proposition~2.21]{Hatcher} and the contractibility of $L$, its augmented
complex is exact. Left translation acts freely on singular simplices, so this complex is a free
$\Z[L^\delta]$-resolution of $\Z$. Sending each simplex to its ordered list of vertices gives an equivariant,
augmentation-preserving chain map to the homogeneous bar resolution and hence an isomorphism in group cohomology.
The pullback of $c$ vanishes on every generator of the small chain complex, so $[c]=0$.
\end{proof}

By Lemma~\ref{h3:geo:diffeo-contractible}, $\cG$ is contractible. Proposition~\ref{apps:localextension} and Lemma~\ref{app:germ} thus show that every real-valued global
homogeneous three-cocycle with the same germ as $\phi_T$ has ordinary class $\xi_{M,3}$. Nosaka's real-valued
cocycle $\Phi_M$ has this germ when the same local fillings $T_k$ are used, since the coefficient in
equation~(9) of that paper is $|\pi_1(\cG)|^3=1$ \cite[Proposition~3.6, \S6.1]{Nosaka}. Hence $[\Phi_M]=\xi_{M,3}$.

Put $D=\Diff_v(M)$, $j:\cG\hookrightarrow D$, and $d=[D:\cG]$. Nosaka's transfer sums over cosets \cite[\S4.2, equation~(14)]{Nosaka}. Apply \eqref{apps:sum-transfer} of Corollary~\ref{apps:hyperbolic-three-transfer} with $u=1$ to
evaluate the transferred class on the integral class $j_*\alpha_{\cG}^{E_{\log}}(1)\in H_3(BD^\delta;\Z)$. This gives
\[
 \big\langle\operatorname{cor}_{\cG}^{D}[\Phi_M],\ j_*\alpha_{\cG}^{E_{\log}}(1)\big\rangle=d\,\totalVol>0 .
\]
In particular, the transferred class $\operatorname{cor}_{\cG}^{D}[\Phi_M]$ is nonzero.
This answers the question posed in \cite[\S6.1]{Nosaka}.


\begin{thebibliography}{99}\small\raggedright
\interlinepenalty=10000
\setlength{\itemsep}{2pt}\setlength{\parsep}{0pt}\setlength{\parskip}{0pt}

\bibitem{McDuffI}
D. McDuff, \emph{Local homology of groups of volume preserving diffeomorphisms. I}, Ann. Sci. \'Ecole Norm. Sup. (4) \textbf{15} (1982), 609--648. \href{https://doi.org/10.24033/asens.1438}{DOI: 10.24033/asens.1438}.

\bibitem{McDuffIII}
D. McDuff, \emph{Local homology of groups of volume-preserving diffeomorphisms. III}, Ann. Sci. \'Ecole Norm. Sup. (4) \textbf{16} (1983), 529--540. \href{https://doi.org/10.24033/asens.1459}{DOI: 10.24033/asens.1459}.

\bibitem{McDuffCan}
D. McDuff, \emph{Some canonical cohomology classes on groups of volume preserving diffeomorphisms}, Trans. Amer. Math. Soc. \textbf{275} (1983), 345--356. \href{https://doi.org/10.1090/S0002-9947-1983-0678355-7}{DOI: 10.1090/S0002-9947-1983-0678355-7}.

\bibitem{Moskowitz}
I. Moskowitz, \emph{Nonvanishing local cohomology classes}, Trans. Amer. Math. Soc. \textbf{286} (1984), 831--837. \href{https://doi.org/10.1090/S0002-9947-1984-0760990-X}{DOI: 10.1090/S0002-9947-1984-0760990-X}.

\bibitem{Hurder}
S. Hurder, \emph{Global invariants for measured foliations}, Trans. Amer. Math. Soc. \textbf{280} (1983), 367--391. \href{https://doi.org/10.1090/S0002-9947-1983-0712266-3}{DOI: 10.1090/S0002-9947-1983-0712266-3}. \href{https://homepages.math.uic.edu/~hurder/papers/11manuscript.pdf}{Author's copy}.

\bibitem{MayPonto}
J. P. May and K. Ponto, \emph{More Concise Algebraic Topology: Localization, Completion, and Model Categories}, Chicago Lectures in Mathematics, University of Chicago Press, Chicago, 2012. \href{https://www.maths.ed.ac.uk/~v1ranick/papers/mayponto.pdf}{Public manuscript}.

\bibitem{BCM}
R. R. Bruner, M. J. Catanzaro and J. P. May, \emph{Characteristic classes}, lecture notes. \href{https://math.uchicago.edu/~may/REU2020/charclasses.pdf}{Public manuscript}.

\bibitem{KKM}
{\latinfont J. K\k{e}dra}, D. Kotschick and S. Morita, \emph{Crossed flux homomorphisms and vanishing theorems for flux groups}, Geom. Funct. Anal. \textbf{16} (2006), 1246--1273. \href{https://arxiv.org/abs/math/0503230}{arXiv:math/0503230}.

\bibitem{Nariman}
S. Nariman, \emph{On the moduli space of flat symplectic surface bundles}, J. Differential Geom. \textbf{116} (2020), 349--391. \href{https://arxiv.org/abs/1612.00043v2}{arXiv:1612.00043v2}.

\bibitem{Moller}
J. M. M{\o}ller, \emph{Nilpotent spaces of sections}, Trans. Amer. Math. Soc. \textbf{303} (1987), 733--741. \href{https://doi.org/10.2307/2000693}{DOI: 10.2307/2000693}. \href{https://web.math.ku.dk/~moller/reprints/nilpotent.pdf}{Author's copy}.

\bibitem{Moore}
J. T. Moore, \emph{Math 7870: Descriptive Set Theory}, lecture notes. \href{https://pi.math.cornell.edu/~justin/Limited2Cornell/DST_lecture_notes.pdf}{Public manuscript}.

\bibitem{Dupont}
J. L. Dupont, \emph{Curvature and Characteristic Classes}, Lecture Notes in Mathematics \textbf{640}, Springer-Verlag, Berlin--New York, 1978. \href{https://webhomes.maths.ed.ac.uk/~v1ranick/papers/dupont2.pdf}{Public copy}.

\bibitem{BM}
M. Brandenbursky and M. Marcinkowski, \emph{Volume and Euler classes in bounded cohomology of transformation groups}, Glasg. Math. J. \textbf{67} (2025), 34--49. \href{https://doi.org/10.1017/S0017089524000223}{DOI: 10.1017/S0017089524000223}. \href{https://arxiv.org/abs/2308.15974}{arXiv:2308.15974}.

\bibitem{Reznikov}
A. Reznikov, \emph{Continuous cohomology of the group of volume-preserving and symplectic diffeomorphisms, measurable transfer and higher asymptotic cycles}, Selecta Math. (N.S.) \textbf{5} (1999), 181--198. \href{https://doi.org/10.1007/s000290050046}{DOI: 10.1007/s000290050046}. \href{https://arxiv.org/abs/dg-ga/9709010}{arXiv:dg-ga/9709010}.

\bibitem{Nosaka}
T. Nosaka, \emph{Chern--Simons-type 3-cocycles and $\mathbb Q/\mathbb Z$-torsion in the third homology of discrete diffeomorphism groups}, \href{https://arxiv.org/abs/2601.07230v3}{arXiv:2601.07230v3} (8 August 2026).

\bibitem{Gabai}
D. Gabai, \emph{The Smale conjecture for hyperbolic 3-manifolds: $\operatorname{Isom}(M^3)\simeq\operatorname{Diff}(M^3)$}, J. Differential Geom. \textbf{58} (2001), 113--149. \href{https://doi.org/10.4310/jdg/1090348284}{DOI: 10.4310/jdg/1090348284}.

\bibitem{MEF}
J. E. Marsden, D. G. Ebin and A. E. Fischer, \emph{Diffeomorphism groups, hydrodynamics and relativity},
in J. Vanstone (ed.), \emph{Proc. 13th Biennial Seminar of the Canadian Mathematical Congress}, 1972, 135--279.
\href{https://www.cds.caltech.edu/~marsden/bib/1972/02-MaEbFi1972/MaEbFi1972.pdf}{Public copy}.

\bibitem{Tshishiku}
B. Tshishiku, \emph{Pontryagin classes of locally symmetric manifolds}, Algebr. Geom. Topol. \textbf{15} (2015), 2707--2754. \href{https://doi.org/10.2140/agt.2015.15.2707}{DOI: 10.2140/agt.2015.15.2707}.

\bibitem{Matsushima}
Y. Matsushima, \emph{On Betti numbers of compact, locally symmetric Riemannian manifolds}, Osaka Math. J. \textbf{14} (1962), 1--20. \href{https://doi.org/10.18910/5711}{DOI: 10.18910/5711}.

\bibitem{Franke}
J. Franke, \emph{A topological model for some summand of the Eisenstein cohomology of congruence subgroups}, in W. T. Gan, S. S. Kudla and Y. Tschinkel (eds.), \emph{Eisenstein Series and Applications}, Progress in Mathematics \textbf{258}, Birkh\"auser Boston, Boston, 2008, 27--85. \href{https://doi.org/10.1007/978-0-8176-4639-4_2}{DOI: 10.1007/978-0-8176-4639-4\_2}. \href{https://math.nyu.edu/~tschinke/princeton/papers/.aim05/submitted/franke.pdf}{Author's manuscript}.

\bibitem{BorelDensity}
A. Borel, \emph{Density properties for certain subgroups of semi-simple groups without compact components}, Ann. of Math. (2) \textbf{72} (1960), 179--188. \href{https://doi.org/10.2307/1970150}{DOI: 10.2307/1970150}.

\bibitem{KellerhalsZehrt}
R. Kellerhals and T. Zehrt, \emph{The Gauss--Bonnet formula for hyperbolic manifolds of finite volume}, Geom. Dedicata \textbf{84} (2001), 49--62. \href{https://doi.org/10.1023/A:1010352829732}{DOI: 10.1023/A:1010352829732}. \href{https://homeweb.unifr.ch/kellerha/pub/GaussBonnet.pdf}{Author's manuscript}.

\bibitem{Hatcher}
A. Hatcher, \emph{Algebraic Topology}, Cambridge University Press, Cambridge, 2002.
\href{https://pi.math.cornell.edu/~hatcher/AT/AT.pdf}{Author's public version}.
\href{https://pi.math.cornell.edu/~hatcher/AT/ATch5.pdf}{Supplementary Chapter 5}.

\bibitem{JHCWhitehead}
J. H. C. Whitehead, \emph{A certain exact sequence}, Ann. of Math. (2) \textbf{52} (1950), 51--110.
\href{https://doi.org/10.2307/1969511}{DOI: 10.2307/1969511}.

\bibitem{Waldschmidt1978}
M. Waldschmidt, \emph{P\'olya's theorem by Schneider's method},
Acta Math. Acad. Sci. Hungar. \textbf{31} (1978), 21--25.
\href{https://webusers.imj-prg.fr/~michel.waldschmidt/articles/pdf/ActaAcadSciHungar31-1978.pdf}{Author's copy of the published article}.


\bibitem{FrigerioBC}
R. Frigerio, \emph{Bounded Cohomology of Discrete Groups}, Mathematical Surveys and Monographs \textbf{227}, AMS, 2017. \href{https://doi.org/10.1090/surv/227}{DOI: 10.1090/surv/227}. \href{https://arxiv.org/abs/1610.08339}{arXiv:1610.08339}.

\bibitem{McCulloughSoma}
D. McCullough and T. Soma, \emph{The Smale conjecture for Seifert fibered spaces with hyperbolic base orbifold}, J. Differential Geom. \textbf{93} (2013), 327--353. \href{https://doi.org/10.4310/jdg/1361800869}{DOI: 10.4310/jdg/1361800869}.

\bibitem{BMInd}
M. Brandenbursky and M. Marcinkowski, \emph{Bounded cohomology of transformation groups}, Math. Ann. \textbf{382} (2022), 1181--1197. \href{https://doi.org/10.1007/s00208-021-02266-8}{DOI: 10.1007/s00208-021-02266-8}.

\bibitem{Mostow}
G. D. Mostow, \emph{Strong Rigidity of Locally Symmetric Spaces}, Annals of Mathematics Studies \textbf{78}, Princeton University Press, 1973.
\bibitem{BFJ}
M. Bustamante, F. T. Farrell and Y. Jiang,
\emph{Rigidity and characteristic classes of smooth bundles with nonpositively curved fibers},
J. Topol. \textbf{9} (2016), 934--956.
\href{https://doi.org/10.1112/jtopol/jtw015}{DOI: 10.1112/jtopol/jtw015}.

\bibitem{ThomRealization}
R. Thom, \emph{Quelques propriétés globales des variétés différentiables},
Comment. Math. Helv. \textbf{28} (1954), 17--86.
\href{https://doi.org/10.1007/BF02566923}{DOI: 10.1007/BF02566923}.

\bibitem{HTV}
S. Haller, J. Teichmann and C. Vizman, \emph{Totally geodesic subgroups of diffeomorphisms},
J. Geom. Phys. \textbf{42} (2002), 342--354.
\href{https://doi.org/10.1016/S0393-0440(01)00096-1}{DOI: 10.1016/S0393-0440(01)00096-1}.

\bibitem{BerglundMadsen}
A. Berglund and I. Madsen, \emph{Rational homotopy theory of automorphisms of manifolds}, Acta Math. \textbf{224} (2020), 67--185. \href{https://doi.org/10.4310/ACTA.2020.v224.n1.a2}{DOI: 10.4310/ACTA.2020.v224.n1.a2}.

\bibitem{GromovVolume}
M. Gromov, \emph{Volume and bounded cohomology}, Publ. Math. Inst. Hautes \'Etudes Sci. \textbf{56} (1982), 5--99.
\href{https://www.numdam.org/item/PMIHES_1982__56__5_0/}{Numdam}.

\end{thebibliography}
\end{document}